\documentclass[reqno,11pt,twoside]{article}
\usepackage[hmargin=1in, vmargin=1.25in, a4paper, centering]{geometry}
\usepackage{amsmath, amsthm,amssymb,enumerate}
\usepackage{fourier}
\usepackage{ebgaramond}
\usepackage{graphicx}
\usepackage{xcolor}
\usepackage[font=small, labelfont=bf, labelsep=colon, justification=centerlast, margin=1in]{caption}
\usepackage{etoolbox, apptools}
\usepackage{longtable,booktabs,array,multirow}
\usepackage{algorithm}
\usepackage{algpseudocode}
\usepackage{marvosym}
\usepackage[section]{placeins}
\usepackage{titletoc}
\numberwithin{equation}{section}
\theoremstyle{plain}
\newtheorem{theorem}{Theorem}[section]
\newtheorem{lemma}[theorem]{Lemma}
\newtheorem{corollary}[theorem]{Corollary}
\theoremstyle{definition}
\newtheorem{ex}[theorem]{Example}
\theoremstyle{remark}
\newtheorem{remark}[theorem]{Remark}
\newcommand{\addQEDstyle}[2]{\AtBeginEnvironment{#1}{\pushQED{\qed}\renewcommand{\qedsymbol}{#2}}\AtEndEnvironment{#1}{\popQED}}
\addQEDstyle{remark}{$\blacktriangleleft$}
\addQEDstyle{ex}{$\blacktriangleleft$}

\newcommand{\dr}{\mathrm{d}}
\newcommand{\ir}{\mathrm{i}}

\renewcommand{\epsilon}{\varepsilon}

\newcommand{\Lfloor}{\left\lfloor}
\newcommand{\Rfloor}{\right\rfloor}
\newcommand{\entire}[1]{\Lfloor #1 \Rfloor}

\renewcommand{\phi}{\varphi}

\usepackage{accents}
\DeclareMathAccent{\wtilde}{\mathord}{largesymbols}{"65}
\DeclareRobustCommand{\utilde}[1]{\underaccent{\wtilde}{#1}}

\renewcommand\footnotemark{}
\usepackage[nobottomtitles,pagestyles]{titlesec}
\titleformat{\section}
{\normalfont\large\bfseries}
{\filcenter\IfAppendix{Appendix }{\S}\thesection.}{1ex}{\filcenter}
\titleformat{\subsection}
{\normalfont\bfseries}
{\filcenter\S\thesubsection.}{1ex}{\filcenter}
\widenhead*{0.5in}{0.5in}
\renewpagestyle{plain}[\bfseries]{\footrule\setfoot{}{\thepage}{}}
\newpagestyle{mystyle}[\bfseries]{
\headrule
\sethead[\thepage][][Uniform enclosures for the phase and zeros of Bessel functions and their derivatives]{N. Filonov, M. Levitin, I. Polterovich, and D. A. Sher}{}{\thepage}
}
\usepackage[colorlinks,allcolors=blue,pagebackref]{hyperref}
\usepackage{pifont}
\renewcommand*{\backrefalt}[4]{%
\ifcase #1 %
No citations%
\or
\ding{43}~p.~#2%
\else
\ding{43}~pp.~#2%
\fi}
\newcommand{\mydoi}[1]{\href{https://doi.org/#1}{doi: #1}}

\usepackage{marginnote}

\begin{document}
\pagestyle{mystyle}
\thispagestyle{plain}
\title{%
Uniform enclosures for the phase and zeros of Bessel functions and their derivatives%
\footnote{The accompanying \texttt{Mathematica} script and its printout are available for download at \url{https://michaellevitin.net/bessels.html}.}
}
\author{
Nikolay Filonov
\thanks{%
\textbf{N. F.: }St. Petersburg Department
of Steklov Institute of Mathematics of RAS,
Fontanka 27, 191023, St.Petersburg, Russia;
St. Petersburg State University,
University emb. 7/9,       
199034, St.Petersburg, Russia; 
\href{mailto:filonov@pdmi.ras.ru}{\nolinkurl{filonov@pdmi.ras.ru}}%
}
\and
Michael Levitin\hspace{-3ex}
\thanks{%
\textbf{M. L.: }Department of Mathematics and Statistics, University of Reading, 
Pepper Lane, Whiteknights, Reading RG6 6AX, UK;
\href{mailto:M.Levitin@reading.ac.uk}{\nolinkurl{M.Levitin@reading.ac.uk}}; \url{https://www.michaellevitin.net}%
}
\and 
Iosif Polterovich
\thanks{%
\textbf{I. P.: }D\'e\-par\-te\-ment de math\'ematiques et de statistique, Univer\-sit\'e de Mont\-r\'eal, 
CP 6128 succ Centre-Ville, Mont\-r\'eal QC  H3C 3J7, Canada;
\href{mailto:iossif@dms.umontreal.ca}{\nolinkurl{iossif@dms.umontreal.ca}}; \url{https://www.dms.umontreal.ca/\~iossif}%
}
\and
David A. Sher
\thanks{%
\textbf{D. A. S.:  }Department of Mathematical Sciences, DePaul University, 2320 N. Kenmore Ave, 60614, Chicago, IL, USA;
\href{mailto:dsher@depaul.edu}{\nolinkurl{dsher@depaul.edu}}
}
}
\date{\small arXiv:2402.06956v3; 1 November 2024\\to appear in SIAM J. Math. Anal. \mydoi{10.1137/24M1642032}}
\maketitle

\begin{abstract} 
We prove explicit uniform two-sided bounds for the phase functions of Bessel functions 
and of their derivatives. As a consequence, we obtain new enclosures for the zeros of Bessel functions and their derivatives in terms of 
inverse values of some elementary functions. These bounds are valid, 
with a few exceptions, for all zeros and all Bessel functions with non-negative indices.  We provide numerical evidence  showing  
that our bounds either improve or closely match the best previously known ones.  
\end{abstract}

{\small\tableofcontents}

\section{Introduction and main results}

\subsection{Setup I: Bessel functions}

Throughout, 
\begin{equation*}
J_\nu(x),\qquad Y_\nu(x)
\end{equation*}
are standard Bessel functions with $\nu\ge 0$, and 
\begin{equation*}
j_{\nu, k},\qquad  y_{\nu, k}
\end{equation*}
are their $k$th largest positive zeros, respectively, where $k\in\mathbb{N}$.  

\begin{figure}[ht]
\centering
\includegraphics{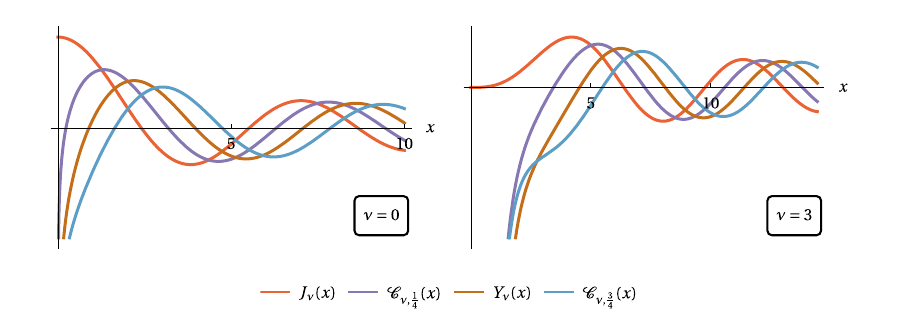}
\caption{Plots of $\mathcal{C}_{\nu, \tau}(x)$.\label{fig:Cnutau}}
\end{figure}

We further define 
\begin{equation*}
\mathcal{C}_{\nu, \tau}(x):=J_\nu(x) \cos(\pi \tau)+ Y_\nu(x)\sin(\pi \tau),
\end{equation*}
where we choose the additional parameter $\tau\in(0, 1]$ for definiteness, see Figure \ref{fig:Cnutau}.
Note that 
\begin{equation*}
J_\nu(x) = - \mathcal{C}_{\nu, 1}(x),\qquad Y_\nu(x) = \mathcal{C}_{\nu, \frac{1}{2}}(x).
\end{equation*}

We also introduce the Bessel \emph{modulus} and \emph{phase functions} $M_\nu, \theta_\nu:(0,+\infty)\to\mathbb{R}$ defined by 
\begin{equation*}
M_\nu(x)=\sqrt{J^2_\nu(x)+Y^2_\nu(x)},\qquad J_\nu(x)+\ir Y_\nu(x)=M_\nu(x)\left(\cos\theta_\nu(x)+\ir \sin\theta_\nu(x)\right),
\end{equation*}
where we choose a continuous branch of $\theta_\nu(x)$ with the initial condition 
\begin{equation*}
\lim_{x\to 0^+}\theta_\nu(x)=-\frac{\pi}{2}.
\end{equation*}
For more details on the modulus and phase functions see \cite[\S10.18]{DLMF} and \cite{Hor}.

The following properties of the Bessel zeros and the modulus and phase functions are standard. We have
\begin{equation*}
\nu<y_{\nu,1}<j_{\nu,1}<y_{\nu,2}<j_{\nu,2}<\dots<y_{\nu,k}<j_{\nu,k}<\dots,
\end{equation*}
see also \cite{PA11} for further interlacing properties. Also,
\begin{equation*}
M_\nu(x)>0,\qquad \theta'_\nu(x)>0\qquad\text{for all }x>0,
\end{equation*}
which means that the inverse function 
\begin{equation*}
\theta_\nu^{-1}: \left[-\frac{\pi}{2}, +\infty\right)\to [0,+\infty)
\end{equation*} 
is well-defined and continuous. We also have, by \cite[(10.18.18)]{DLMF}, the asymptotics
\begin{equation}\label{eq:thetaasympt}
\theta_{\nu}(x) = x-\frac{\pi}{4}\left(2\nu+1\right)+%
\frac{4\nu^2-1}{8x}+\frac{(4\nu^2-1)(4\nu^2-25)}{384x^3}+
O\left(x^{-5}\right)\qquad\text{as }x\to+\infty.
\end{equation}

\begin{figure}[ht]
\centering
\includegraphics{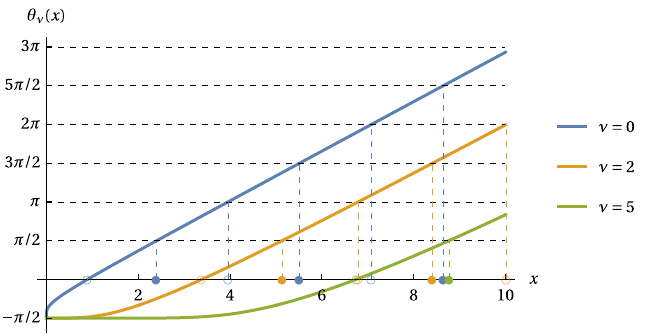}
\caption{Plots of $\theta_\nu(x)$. The filled colour-coded dots on the horizontal axis indicate the positions of zeros $j_{\nu,k}$ and the hollow dots the positions of zeros $y_{\nu,k}$. The phase functions are calculated using the method of \cite{Hor}.\label{fig:theta}}
\end{figure}

We have, since $j_{\nu,k}$ are the positive zeros of $\cos\theta_\nu(x)$, and $y_{\nu,k}$ are the positive zeros of $\sin\theta_\nu(x)$, the relations
\begin{equation}\label{eq:jyk}
j_{\nu,k} = \theta_\nu^{-1}\left(\pi\left(k-\frac{1}{2}\right)\right),\qquad y_{\nu,k} = \theta_\nu^{-1}\left(\pi\left(k-1\right)\right),\qquad k\in\mathbb{N}.
\end{equation}
Some typical plots of functions $\theta_\nu$ are shown in Figure \ref{fig:theta}.

We also have
\begin{equation*}
\mathcal{C}_{\nu, \tau}(x)=M_\nu(x)\cos\left(\theta_\nu(x)-\pi \tau\right).
\end{equation*}
The $k$th positive zero of this function, which we will denote by $c_{\nu,\tau, k}$, is given by 
\begin{equation}\label{eq:cnutauk}
c_{\nu,\tau, k} = \theta_\nu^{-1}\left(\pi \left(\tau+k-\frac{3}{2}\right)\right).
\end{equation}

\begin{figure}[ht]
\centering
\includegraphics{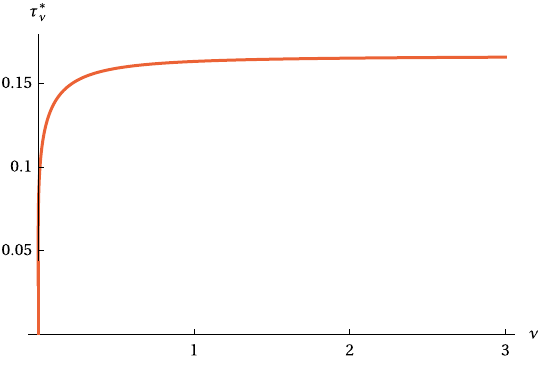}
\caption{Plot of $\tau^*_\nu$ against $\nu$.\label{fig:taustar}}
\end{figure}

\begin{remark} We know that the first zero of either $J_\nu(x)$ or $Y_\nu(x)$ is always bigger that $\nu$. This is not necessarily the case for  $c_{\nu,\tau, 1}$ with $\tau\in \left(0,\frac{1}{2}\right)$. Indeed, we have 
\begin{equation*}
c_{\nu,\tau, 1}
= \theta_\nu^{-1}\left(\pi \left(\tau-\frac{1}{2}\right)\right)\le \nu
\end{equation*}
whenever 
\begin{equation*}
0<\tau\le \tau^*_\nu := \frac{1}{\pi}\theta_\nu(\nu)+\frac{1}{2}. 
\end{equation*} 
Since $\theta_\nu(\nu)\in\left[-\frac{\pi}{2},0\right)$, we always have $\tau^*_\nu\in \left[0,\frac{1}{2}\right)$, see Figure \ref{fig:taustar}.
\end{remark}

For further use, we additionally define the \emph{counting function} of Bessel zeros,
\begin{equation}\label{eq:NJ}
\mathcal{N}_{J_\nu}(\lambda):=\#\left\{k\in\mathbb{N}: j_{\nu,k}\le \lambda\right\}.
\end{equation}

\subsection{Setup II: derivatives of Bessel functions}

\begin{figure}[ht]
\centering
\includegraphics{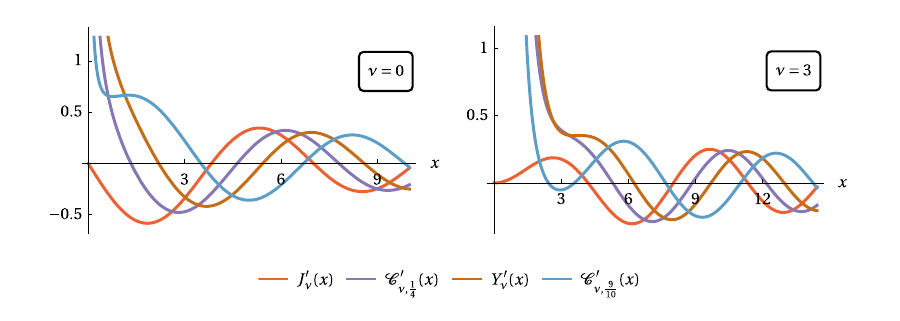}
\caption{Plots of $\mathcal{C}'_{\nu, \tau}(x)$.\label{fig:Cnutauprime}}
\end{figure}

Consider, for $\nu\ge 0$,  the derivatives of the Bessel functions 
\begin{equation*}
J'_\nu(x),\qquad Y'_\nu(x).
\end{equation*}
The numbers
\begin{equation*}
j'_{\nu, k},\qquad  y'_{\nu, k}
\end{equation*}
are their $k$th largest positive zeros except when $\nu=0$, in which case we set
\begin{equation*}
 j'_{0, 1}:=0.
\end{equation*}
We also consider the linear combinations
\begin{equation*}
\mathcal{C}'_{\nu, \tau}(x):=J'_\nu(x) \cos(\pi \tau)+ Y'_\nu(x)\sin(\pi \tau),
\end{equation*}
this time with $\tau\in[0, 1)$, see Figure \ref{fig:Cnutauprime}.

As for the Bessel functions themselves, we also introduce the \emph{modulus} and \emph{phase functions} of the Bessel derivatives, $N_\nu, \phi_\nu:(0,+\infty)\to\mathbb{R}$ defined by 
\begin{equation*}
N_\nu(x)=\sqrt{\left(J'_\nu(x)\right)^2+\left(Y'_\nu(x)\right)^2},\qquad J'_\nu(x)+\ir Y'_\nu(x)=N_\nu(x)\left(\cos\phi_\nu(x)+\ir \sin\phi_\nu(x)\right),
\end{equation*}
where we choose a continuous branch of $\phi_\nu(x)$ with the initial condition 
\begin{equation*}
\lim_{x\to 0^+}\phi_\nu(x)=\frac{\pi}{2},
\end{equation*}
see \cite[\S10.18]{DLMF} and \cite{Hor}.

The following properties of the Bessel derivatives zeros and the modulus and phase functions are standard:
\begin{equation*}
\nu\le j'_{\nu,1}<y'_{\nu,1}<j'_{\nu,2}<y'_{\nu,2}<\dots<j'_{\nu,k}<y'_{\nu,k}<\dots
\end{equation*}
(with equality only for  $\nu=0$);
\begin{equation*}
N_\nu(x)>0\qquad\text{for all }x>0.
\end{equation*}
The phase function $\phi_\nu$ is not monotone (except for $\nu=0$) but has a single minimum at $x=\nu$, with $\phi_\nu(\nu)\in \left(0,\frac{\pi}{2}\right)$ for $\nu>0$.
This means that its inverse function may only be defined on $\left[\phi_\nu(\nu), +\infty\right)$; it is however easier to further restrict its domain and treat it as the function 
\begin{equation*}
\phi_\nu^{-1}: \left[\frac{\pi}{2}, +\infty\right)\to \left[j'_{\nu,1},+\infty\right).
\end{equation*} 
Additionally, by \cite[(10.18.21)]{DLMF},
\begin{equation}\label{eq:phiasympt}
\phi_\nu(x)= x-\frac{\pi}{4}\left(2\nu+1\right)+%
\frac{4\nu^2+3}{8x}+\frac{16\nu^4+184\nu^2-63}{384x^3}+O\left(x^{-5}\right)\qquad\text{as }x\to+\infty.
\end{equation}

We have
\begin{equation}\label{eq:jynupk}
j'_{\nu,k} =\phi_\nu^{-1}\left(\pi\left(k-\frac{1}{2}\right)\right),\qquad y'_{\nu,k} =\phi_\nu^{-1}\left(\pi k\right),\qquad k\in\mathbb{N}.
\end{equation}
Further, let
\begin{equation}\label{eq:ctnupk}
c'_{\nu,\tau, k}:=
\phi_\nu^{-1}\left(\pi \left(\tau+k-\frac{1}{2}\right)\right), \qquad k\in\mathbb{N},
\end{equation}
be the $k$th largest zero of 
\begin{equation*}
\mathcal{C}'_{\nu, \tau}(x)=N_\nu(x)\cos\left(\phi_\nu(x)-\pi \tau\right).
\end{equation*}
in the interval  $\left[j'_{\nu,1}, +\infty\right)$.

\begin{ex} Our choice of enumeration of $c'_{\nu,\tau, k}$ means that we may disregard some positive zeros of $\mathcal{C}'_{\nu, t}$ below $j'_{\nu,1}$. For example, the right-hand side plot in Figure \ref{fig:Cnutauprime} shows that there are two zeros of $\mathcal{C}'_{3, \frac{9}{10}}$ below the first zero of $J'_3$, which are therefore excluded. 
\end{ex}

\begin{figure}[ht]
\centering
\includegraphics{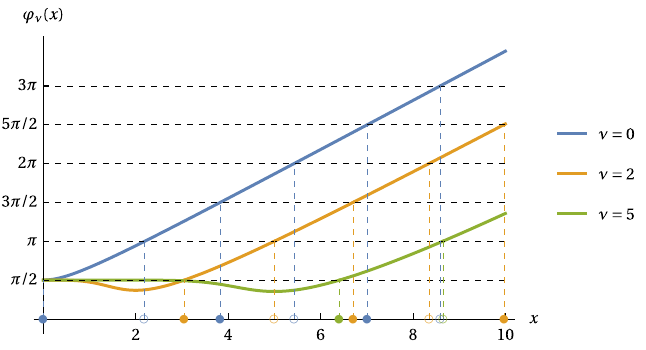}
\caption{Plots of $\phi_\nu(x)$. The filled colour-coded dots on the horizontal axis indicate the positions of zeros $j'_{\nu,k}$ and the hollow dots the positions of zeros $y'_{\nu,k}$. The phase functions are calculated using the method of \cite{Hor}.\label{fig:phi}}
\end{figure}

Some typical plots of functions $\phi_\nu$
are shown in Figure \ref{fig:phi}.

Similarly to \eqref{eq:NJ}, we also define the \emph{counting function} of zeros of Bessel derivatives,
\begin{equation}\label{eq:NJprime}
\mathcal{N}_{J'_\nu}(\lambda):=\#\left\{k\in\mathbb{N}: j'_{\nu,k}\le \lambda\right\}.
\end{equation}

\subsection{Our philosophy and some history}

The main purpose of our paper is to present new uniform (in $\nu$) two-sided bounds for the phase functions $\theta_\nu(x)$ and $\phi_\nu(x)$ (Theorems \ref{thm:thetabounds} and \ref{thm:phibounds}), and, as a consequence, new bounds for the zeros of Bessel functions and their derivatives (Corollaries \ref{cor:jbounds} and \ref{cor:jprimebounds}) valid, with minimal exceptions, for all $\nu$ and $k$. Our bounds for the phase functions are expressed in terms of explicit elementary functions  defined in \S\ref{subs:auxI}  and \S\ref{subs:auxII}, and the bounds for the zeros in terms of inverses of these functions. Additionally, due to the elementary nature of our bounds for phase functions, we derive simple bounds for the counting functions of zeros of Bessel functions and their derivatives (Corollaries \ref{cor:NJbounds} and \ref{cor:NJprimebounds}), some of which have already been used to great effect in spectral-geometric
problems, see \cite{FLPS} and  \cite{FLPS-AB}.

The common feature of our bounds and pre-existing bounds on zeros of Bessel functions and of their derivatives is that they are all based, in one form or another, on asymptotics of Bessel functions for large argument and/or order. We recall that there are, roughly speaking, three relevant asymptotic regimes, see \cite[\S\S10.17, 10.19--20]{DLMF}:
\begin{enumerate}
\item[(R$_\text{i}$)] the large-argument region, in which $x$ approaches infinity and $\nu$ is fixed (or at least bounded);
\item[(R$_\text{ii}$)] the transition region, in which $\nu$ approaches infinity with $x\sim \nu+a\nu^{1/3}$ for some $a\in\mathbb R$. In this region the Bessel functions exhibit Airy-type behaviour;
\item[(R$_\text{iii}$)] the Debye expansion region, in which $\nu$ approaches infinity with $x\sim c\nu$ for some $c>1$.
\end{enumerate}   
Among the extensive historical literature on the topic we mention the work of  Langer \cite{Lan} in the Debye region (R$_\text{iii}$) and of  Chester, Friedman, and Ursell \cite{CFU} in the transition region (R$_\text{ii}$). For a modern perspective on the interplay between these regimes, see \cite{Sh}.

The general previously used approach to bounding Bessel zeros is essentially based on using the asymptotics of Bessel zeros in either of the regimes (R$_\text{i}$) or (R$_\text{ii}$), 
and then applying a form of the Sturm comparison theorem.  Some existing bounds for the zeros of Bessel functions obtained using the expansions in the large argument regime (R$_\text{i}$) are due to Hethcote \cite{Heth} and Elbert--Laforgia \cite{EL00}, see also 
\cite{GG00, GG07, F18, Nem21} for related results. The explicit expressions of these bounds, which we use for benchmarking purposes, are collected in \S\ref{sec:benchmark}. Note that the lower bounds for Bessel zeros obtained in this way are only available for very low values of $\nu$. Alternative bounds, based on the asymptotics of $j_{\nu,k}$ in the transitional regime (R$_\text{ii}$), are given by Qu and Wong \cite{QW}, and are also presented in \S\ref{sec:benchmark}. There is significantly less information in the literature on uniform bounds for the zeros $j'_{\nu,k}$ of derivatives of Bessel functions: we failed to find any bounds derived from the asymptotics in the large argument regime (R$_\text{i}$). The upper bound for $j'_{\nu,k}$ derived from the asymptotics in the transitional regime (R$_\text{ii}$)  from \cite{EL97}, \cite[\S1.7]{E01} is listed in \S\ref{sec:benchmark}.

Unlike the previous approaches, we choose to work directly with the phase functions $\phi_\nu(x)$ and $\theta_\nu(x)$ instead, using their non-oscillatory character and their known asymptotics in  the Debye expansion region (R$_\text{iii}$). It has been known that these asymptotics give, in practice, good approximations of Bessel functions and their derivatives, and can be used for effective calculation of corresponding zeros with high precision,  see \cite{Hor}  as well as various generalisations in the important series of papers by Bremer, Rokhlin, et al.\ \cite{HBRV15, HBR15, BR16, B17, B19}.  Using techniques of Olver, Horsley  \cite{Hor} has recently deduced asymptotics for the phase functions $\theta_{\nu}(x)$ and $\phi_\nu(x)$ in the Debye region (R$_\text{iii}$). An extremely helpful feature of this problem is that, as it turns out,  the expansions in the Debye region also give effective expansions in the large-argument region, and only weaken when one proceeds deep into the transition region. Further, it turns out that truncating these expansions leads to lower and upper bounds for phase functions: the error terms are sign-definite for a broad range of values. We prove this 
using a consequence of the Sturm comparison theorem  (Theorem \ref{thm:sturm}) combined with the non-oscillatory nature of the phase functions. The more detailed comparison of our bounds with existing ones and their effectiveness will be addressed in \S\ref{sec:benchmark}.

The rest of the paper is organised as follows. In \S\S\ref{subs:auxI}--\ref{subs:mainII}, we define the necessary auxiliary functions, study their properties, and state our main results, illustrating them graphically. All the proofs, together with additional Sturm--Liouville theory required, are collected in \S\ref{sec:proofs}.
In \S\ref{sec:ultra}, we extend our techniques to obtain the bounds on the zeros of derivatives of \emph{ultraspherical} Bessel functions. As we have mentioned, \S\ref{sec:benchmark} discusses the comparison of our results with older ones, based on numerical evidence collected in Appendix \ref{sec:numerical}. In order not to overload the main text, some complicated intermediate expressions are placed in Appendix \ref{sec:potentials}.

We remark  that many of the proofs involve straightforward but rather cumbersome manipulations with explicit functions. Although we have verified them independently by hand, we have relied in some cases on simplifications using \texttt{Mathematica}. In the interest of transparency, the corresponding script, as well as its printout, containing \emph{all} of the analytic and numerical calculations we performed, including plotting routines, is available for inspection and download at

\centerline{\url{https://michaellevitin.net/bessels.html}.}

\subsection{Definitions and properties of the auxiliary functions I}\label{subs:auxI}
We start by defining the following auxiliary functions. Everywhere below $x\ge \nu\ge 0$.

We define
\begin{equation*}
\wtilde{\theta}_\nu(x):=\sqrt{x^2-\nu^2}-\nu\arccos\frac{\nu}{x}-\frac {\pi}{4},
\end{equation*}
and note that 
\begin{equation*}
\wtilde{\theta}_\nu(\nu)=-\frac {\pi}{4},
\end{equation*}
\begin{equation}\label{eq:overthetaasympt}
\wtilde{\theta}_\nu(x)=x-\frac{\pi}{4}  (2 \nu +1)+\frac{\nu ^2}{2x}+O\left(x^{-3}\right)\qquad\text{as }x\to+\infty.
\end{equation}
We also note that 
\begin{equation*}
\wtilde{\theta}'_\nu(x)=\frac{\sqrt{x^2-\nu^2}}{x}>0\qquad\text{for } x>\nu.
\end{equation*}
Therefore, the function $\wtilde{\theta}_\nu$ is strictly monotone increasing on the interval $(\nu, +\infty)$, and thus the inverse function 
\begin{equation*}
\left(\wtilde{\theta}_\nu\right)^{-1}:\left[-\frac {\pi}{4},+\infty\right)\to\left[\nu,+\infty\right)
\end{equation*}
is well-defined.

Further, we set
\begin{equation*}
\utilde{\theta}_\nu(x):=\wtilde{\theta}_\nu(x)-\frac{3x^2+2\nu^2}{24(x^2-\nu^2)^{3/2}}=\sqrt{x^2-\nu^2}-\nu\arccos\frac{\nu}{x}-\frac {\pi}{4}-\frac{3x^2+2\nu^2}{24(x^2-\nu^2)^{3/2}},
\end{equation*}
and note that 
\begin{equation*}
\utilde{\theta}_\nu(x)\to-\infty\quad\text{as }x\to\nu^+,
\end{equation*}
\begin{equation}\label{eq:underthetaasympt}
\utilde{\theta}_\nu(x)=x - \frac{\pi}{4} (2 \nu +1)+\frac{4 \nu ^2-1}{8 x}
+\frac{\nu ^2 \left(2 \nu ^2-13\right)}{48 x^3}+O\left(x^{-5}\right)\qquad\text{as }x\to+\infty.
\end{equation}
We also have 
\begin{equation*}
{\utilde{\theta}}'_\nu(x)=\frac{8 x^6+\left(1-24 \nu ^2\right) x^4+4 \left(6 \nu ^4+\nu ^2\right) x^2 - 8 \nu ^6}{8 x (x^2-\nu^2)^{5/2}},
\end{equation*}
and therefore 
\begin{equation*}
{\utilde{\theta}}'_\nu\left(\sqrt{\nu^2+\chi}\right)=\frac{8 \chi^3+\chi ^2+6 \chi \nu ^2+5 \nu ^4}{8 \chi^{5/2} \sqrt{\nu ^2+\chi}}>0\qquad\text{for }\chi>0.
\end{equation*}
Hence the function $\utilde{\theta}_\nu$ is strictly monotone increasing on the interval $(\nu, +\infty)$, and the inverse function 
\begin{equation*}
\left(\utilde{\theta}_\nu\right)^{-1}:\mathbb{R}\to\left(\nu,+\infty\right)
\end{equation*}
is well-defined.

\begin{remark}\label{rem:asDebye} We should explain the origin of the functions $\wtilde{\theta}_\nu(x)$ and $\utilde{\theta}_\nu(x)$: they represent, correspondingly, the two- and three-term asymptotic expansions of $\theta_\nu(x)$ in the Debye regime (R$_\text{iii}$), see \cite[formula (34)]{Hor}, which are obtained using the algorithm of \cite{Olv74}.
\end{remark}

\subsection{Main results I: bounding the phase and zeros of Bessel functions}

Set, for $x>\nu\ge 0$,
\begin{equation*}
\utilde{\utilde{\theta}}_\nu(x):=\max\left\{\utilde{\theta}_\nu(x), -\frac{\pi}{2}\right\}.
\end{equation*}
We have

\begin{theorem}\label{thm:thetabounds} For every $\nu\ge 0$ and every $x>\nu$,
\begin{equation}\label{eq:thetabounds}
\utilde{\utilde{\theta}}_\nu(x)<\theta_\nu(x)<\wtilde{\theta}_\nu(x).
\end{equation}
\end{theorem}

\begin{remark} The upper bound in Theorem \ref{thm:thetabounds} has been already proved, using different techniques, in \cite{Sh} and \cite{FLPS}.
\end{remark}

Recalling \eqref{eq:jyk} and \eqref{eq:cnutauk}, Theorem \ref{thm:thetabounds} immediately implies
\begin{corollary}\label{cor:jbounds} 
For every $\nu\ge 0$, we have the following bounds.
\begin{enumerate}[{\normalfont(i)}]
\item For every $k\in\mathbb{N}$,
\begin{equation*}
\utilde{j}_{\nu,k}:=\left(\wtilde{\theta}_\nu\right)^{-1}\left(\pi\left(k-\frac{1}{2}\right)\right)<j_{\nu,k}<\left(\utilde{\theta}_\nu\right)^{-1}\left(\pi\left(k-\frac{1}{2}\right)\right)=:\wtilde{j}_{\nu,k}.
\end{equation*}
\item For every $k\in\mathbb{N}$,
\begin{equation*}
\utilde{y}_{\nu,k}:=\left(\wtilde{\theta}_\nu\right)^{-1}\left(\pi (k-1) \right)<y_{\nu,k}<\left(\utilde{\theta}_\nu\right)^{-1}\left(\pi\ (k-1)\right)=:\wtilde{y}_{\nu,k}.
\end{equation*}
\item More generally, for all $\tau\in(0,1]$ and all $k\ge 2$, we have
\begin{equation*}
\utilde{c}_{\nu,\tau,k}:=\left(\wtilde{\theta}_\nu\right)^{-1}\left(\pi\left(\tau+k-\frac{3}{2}\right)\right)<c_{\nu,\tau,k}<\left(\utilde{\theta}_\nu\right)^{-1}\left(\pi\left(\tau+k-\frac{3}{2}\right)\right)=:\wtilde{c}_{\nu,\tau,k}.
\end{equation*}
For $k=1$, the upper bound $c_{\nu,\tau,1}<\wtilde{c}_{\nu,\tau,1}$ is also valid for all $\tau\in(0,1]$, and the lower bound $\utilde{c}_{\nu,\tau,1}<c_{\nu,\tau,1}$ is valid for $\tau\in\left(\frac{1}{4},1\right]$.
\end{enumerate}
\end{corollary}
In part (iii), the restriction on the range of $\tau$ for which the lower bound holds with $k=1$ is due to the necessary condition $\pi\left(\tau+\frac{1}{2}\right) > -\frac{\pi}{4}$ for applying the inverse function.

\begin{figure}[ht]
\centering
\includegraphics{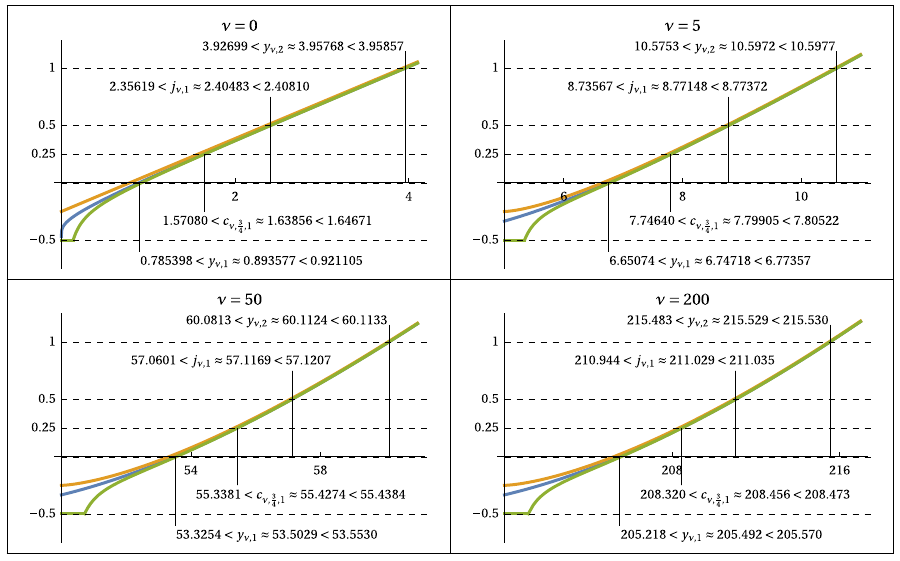}
\caption{All plots show the scaled Bessel phase function $\frac{1}{\pi}\theta_\nu(x)$ (blue), and its scaled bounds $\frac{1}{\pi}\wtilde{\theta}_\nu(x)$ (orange) and  $\frac{1}{\pi}\utilde{\utilde{\theta}}_\nu(x)$ (green). The zeros $c_{\nu, \frac{3}{4}, 1}$ are the first positive zeros of $\mathcal{C}_{\nu, \frac{3}{4}}(x)$, which up to a constant factor coincides with $J_\nu(x)-Y_\nu(x)$.\label{fig:comparisonJ}}
\end{figure}

For an illustration of  the results in Theorem \ref{thm:thetabounds} and Corollary \ref{cor:jbounds}, see Figure \ref{fig:comparisonJ}.

We also immediately deduce the bounds for the counting function \eqref{eq:NJ} of Bessel zeros.

\begin{corollary}\label{cor:NJbounds} 
For any $\lambda>\nu\ge 0$,
\begin{equation*}
\entire{\frac{1}{\pi}\utilde{\utilde{\theta}}_\nu(\lambda)+\frac{1}{2}}
\le  \mathcal{N}_{J_\nu}(\lambda)
\le \entire{\frac{1}{\pi}\wtilde{\theta}_\nu(\lambda)+\frac{1}{2}},
\end{equation*}
where $\entire{\cdot}$ denotes the integer part.
\end{corollary}

\subsection{Definitions and properties of the auxiliary functions II}\label{subs:auxII}
We define, for $x>\nu\ge 0$, the functions
\begin{equation}\label{eq:defunderphi}
\utilde{\phi}_\nu(x):=\wtilde{\theta}_\nu(x)+\frac{\pi}{2}=\sqrt{x^2-\nu^2}-\nu\arccos\frac{\nu}{x}+\frac {\pi}{4},
\end{equation}
and
\begin{equation*}
\wtilde{\phi}_\nu(x):=\utilde{\phi}_\nu(x)+\frac{9x^2-2\nu^2}{24(x^2-\nu^2)^{3/2}}=\sqrt{x^2-\nu^2}-\nu\arccos\frac{\nu}{x}+\frac {\pi}{4}+\frac{9x^2-2\nu^2}{24(x^2-\nu^2)^{3/2}}.
\end{equation*}

As in \S\ref{subs:auxI}, we have 
\begin{equation*}
\utilde{\phi}_\nu(\nu)=\frac{\pi}{4},
\end{equation*}
\begin{equation}\label{eq:underphiasympt}
\utilde{\phi}_\nu(x)=x-\frac{\pi}{4}  (2 \nu - 1)+\frac{\nu ^2}{2x}+O\left(x^{-3}\right)\qquad\text{as }x\to+\infty.
\end{equation}
Also,
\begin{equation*}
{\utilde{\phi}}'_\nu(x)=\wtilde{\theta}'_\nu(x)>0\qquad\text{for all }x>\nu\ge 0,
\end{equation*}
and the inverse function
\begin{equation*}
\left(\utilde{\phi}_\nu\right)^{-1}:\left[\frac{\pi}{4},+\infty\right)\to\left[\nu,+\infty\right)
\end{equation*}
is well-defined. Moreover,
\begin{equation*}
\left(\utilde{\phi}_\nu\right)^{-1}(z)=\left(\wtilde{\theta}_\nu\right)^{-1}\left(z-\frac{\pi}{2}\right).
\end{equation*}

The behaviour of the function $\wtilde{\phi}_\nu(x)$ is more complicated. We have 
\begin{equation*}
\wtilde{\phi}_\nu(x)\to+\infty\quad\text{as }x\to\nu^+,
\end{equation*}
and
\begin{equation}\label{eq:overphiasympt}
\wtilde{\phi}_\nu(x)=x-\frac{\pi}{4} (2 \nu -1)+\frac{4 \nu ^2+3}{8 x}+\frac{\nu ^2 \left(2 \nu ^2+23\right)}{48 x^3}+O\left(x^{-5}\right)\qquad\text{as }x\to+\infty.
\end{equation}
Further,
\begin{equation*}
\wtilde{\phi}'_\nu(x)=\frac{p_\nu(x)}{8 x \left(x^2-\nu ^2\right)^{5/2}},
\end{equation*}
where
\begin{equation}\label{eq:pnu}
p_\nu(x):=8 x^6-3 \left(8 \nu ^2+1\right) x^4+4 \nu ^2 \left(6 \nu ^2-1\right) x^2 - 8 \nu ^6
\end{equation}
is not sign-definite for $x\in[\nu,+\infty)$. Namely, we will prove

\begin{lemma}\label{lem:xstar}
For any $\nu\ge 0$, we have $\wtilde{\phi}'_\nu(x)>0$ for $x\in\left(x^\star_\nu,+\infty\right)$, where $x^\star_\nu$ is the only root in $(\nu,+\infty)$ of \eqref{eq:pnu}. The quantity 
\begin{equation*}
z^\star_\nu:=\wtilde{\phi}_\nu\left(x_\nu^\star\right)
\end{equation*}
is positive and monotone decreasing in $\nu$ for all $\nu\ge 0$.
\end{lemma}

\begin{figure}[ht]
\centering
\includegraphics{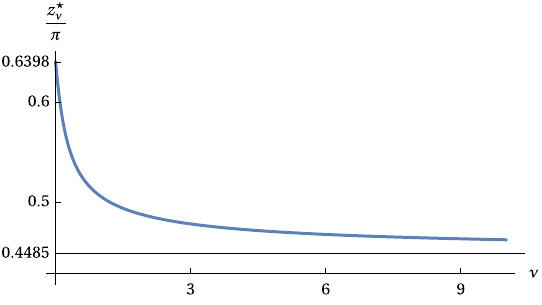}
\caption{The plot of  $\frac{1}{\pi}z^\star_\nu$ against $\nu$, with the horizontal asymptote positioned at an approximate value of $\frac{1}{4}+\frac{1}{\pi}\sqrt{\frac{7}{18}}$.\label{fig:zstar}}
\end{figure}

We present the graph of $\frac{1}{\pi}z^\star_\nu$ as a function of $\nu$ in Figure \ref{fig:zstar}. Note that $x_0^\star=\sqrt{\frac{3}{8}}$ and $z_0^\star=\frac{\pi}{4} + \sqrt{\frac{3}{2}}$.  It may be also shown that
\begin{equation*}
x^\star_\nu=\nu+\frac{\sqrt[3]{7\nu}}{4}+O\left(\nu^{-1/3}\right),\qquad
z^\star_\nu=\frac{\pi}{4}+\sqrt{\frac{7}{18}}+O\left(\nu^{-2/3}\right)\qquad\text{as }\nu\to+\infty.
\end{equation*}

Lemma \ref{lem:xstar} now ensures that the inverse function 
\begin{equation*}
\left(\utilde{\phi}_\nu\right)^{-1}:\left[z_\nu^\star,+\infty\right)\to\left[x_\nu^\star,+\infty\right)
\end{equation*}
is well-defined.

\begin{remark}\label{rem:asDebyeprime} Similarly to Remark \ref{rem:asDebye}, the functions $\utilde{\phi}_\nu(x)$ and $\wtilde{\phi}_\nu(x)$ represent, correspondingly, the two- and three-term asymptotic expansions of $\phi_\nu(x)$ in the Debye regime (R$_\text{iii}$), see \cite[formula (46)]{Hor}.
\end{remark}

\subsection{Main results II: bounding the phase and zeros of derivatives of  Bessel functions}\label{subs:mainII}

Set, for $x>\nu\ge 0$,
\begin{equation*}
\wtilde{\wtilde{\phi}}_\nu(x):=\begin{cases}z^\star_\nu\quad&\text{if }x<x^\star_\nu,\\\wtilde{\phi}_\nu(x)\quad&\text{if }x\ge x^\star_\nu.\end{cases}
\end{equation*}

\begin{theorem}\label{thm:phibounds}  For every $\nu\ge 0$ and every $x>\nu$, we have
\begin{equation}\label{eq:phibounds}
\utilde{\phi}_\nu(x)<\phi_\nu(x)<\wtilde{\wtilde{\phi}}_\nu(x).
\end{equation}
\end{theorem}

\begin{remark} 
The lower bound in Theorem \ref{thm:phibounds} has been already proved, using different techniques, in \cite{FLPS}.
\end{remark}

Recalling \eqref{eq:jynupk} and \eqref{eq:ctnupk}, Theorem \ref{thm:phibounds} immediately implies

\begin{figure}[ht]
\centering
\includegraphics{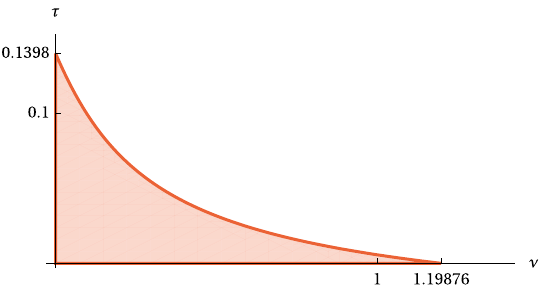}
\caption{The region of the $(\nu,\tau)$-plane in which the lower bound $\utilde{c}'_{\nu,\tau,1}<c'_{\nu,\tau,1}$  for the first zero is \emph{not applicable}.\label{fig:forbid}}
\end{figure}

\begin{corollary}\label{cor:jprimebounds} 
For every $\nu\ge 0$, we have the following bounds.
\begin{enumerate}[{\normalfont(i)}]
\item For every $k\in\{2, 3, \dots\}$,
\begin{equation*}
\utilde{j}'_{\nu,k}:=\left(\wtilde{\phi}_\nu\right)^{-1}\left(\pi\left(k-\frac{1}{2}\right)\right)<j'_{\nu,k}<\left(\utilde{\phi}_\nu\right)^{-1}\left(\pi\left(k-\frac{1}{2}\right)\right)=:\wtilde{j}'_{\nu,k}.
\end{equation*}
For $k=1$, the upper bound $j'_{\nu,1}<\wtilde{j}'_{\nu,1}$ is also valid for all $\nu\ge 0$, and the lower bound $\utilde{j}'_{\nu,1}<j'_{\nu,1}$ is valid whenever
\begin{equation*}
z^\star_\nu\le \frac{\pi}{2},
\end{equation*}
that is, roughly, for $\nu\gtrsim 1.19876$.
\item For every $k\in\mathbb{N}$,
\begin{equation*}
\utilde{y}'_{\nu,k}:=\left(\wtilde{\phi}_\nu\right)^{-1}(\pi k)<y'_{\nu,k}<\left(\utilde{\phi}_\nu\right)^{-1}(\pi k)=:\wtilde{y}'_{\nu,k}.
\end{equation*}
\item More generally, for all $\tau\in[0,1)$ and all $k\ge 2$, we have
\begin{equation*}
\utilde{c}'_{\nu,\tau,k}:=\left(\wtilde{\phi}_\nu\right)^{-1}\left(\pi\left(\tau+k-\frac{1}{2}\right)\right)<c'_{\nu,\tau,k}<\left(\utilde{\phi}_\nu\right)^{-1}\left(\pi\left(\tau+k-\frac{1}{2}\right)\right)=:\wtilde{c}'_{\nu,\tau,k}.
\end{equation*}
For $k=1$, the upper bound $c'_{\nu,\tau,1}<\wtilde{c}'_{\nu,\tau,1}$ is also valid for all $\tau\in[0,1)$. The lower bound $\utilde{c}'_{\nu,\tau,1}<c'_{\nu,\tau,1}$  is valid for $\nu$ such that
\begin{equation*}
z^\star_\nu\le \pi\left(\tau+\frac{1}{2}\right);
\end{equation*}
in particular, if $\tau\ge \frac{1}{\pi} z^\star_0-\frac{1}{2}=\frac{1}{\pi}\sqrt{\frac{3}{2}}-\frac{1}{4}\approx 0.139848$, then it is valid for all $\nu\ge 0$, see Figure \ref{fig:forbid}.
\end{enumerate}
\end{corollary}

In parts (i) and (iii), the restriction on the range of $\nu$ and $\tau$ for which the lower bound holds with $k=1$ is due to the fact that the argument of the inverse function $\left(\wtilde{\phi}_\nu\right)^{-1}$ should exceed $z^\star_\nu$.

\begin{figure}[ht]
\centering
\includegraphics{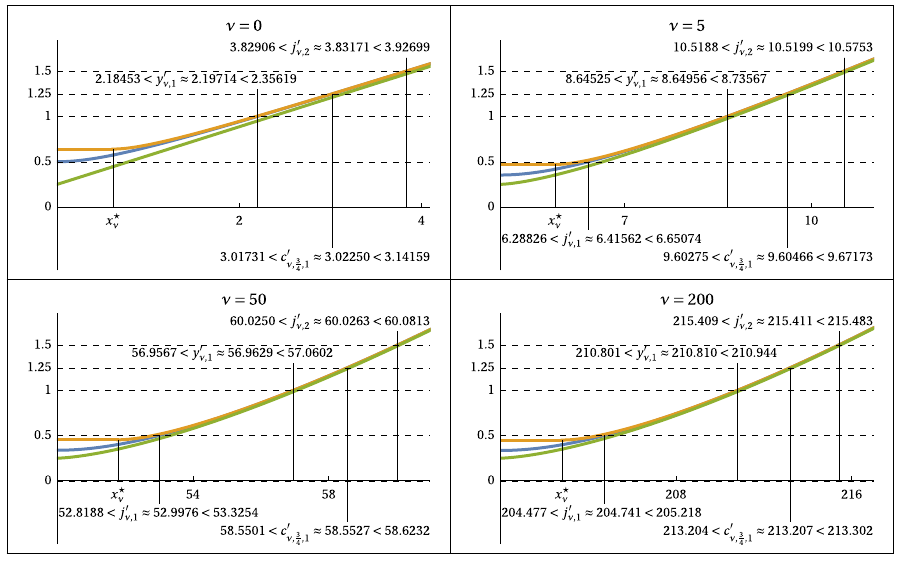}
\caption{All plots show the scaled Bessel derivative phase function $\frac{1}{\pi}\phi_\nu(x)$ (blue), and its scaled bounds $\frac{1}{\pi}\wtilde{\wtilde{\phi}}_\nu(x)$ (orange) and  $\utilde{\phi}_\nu(x)$ (green). The zeros $c'_{\nu, \frac{3}{4}, 1}$ are the first zeros above $j'_{\nu,1}$ of $\mathcal{C}'_{\nu, \frac{3}{4}}(x)$, which up to a constant factor coincides with $J'_\nu(x)-Y'_\nu(x)$.\label{fig:comparisonJp}}
\end{figure}

For an illustration of  the results in Theorem  \ref{thm:phibounds} and Corollary \ref{cor:jprimebounds}, see Figure \ref{fig:comparisonJp}.

We also immediately deduce the bounds for the counting function \eqref{eq:NJprime} of zeros of Bessel derivatives.

\begin{corollary}\label{cor:NJprimebounds} 
For any $\lambda>\nu\ge 0$,
\begin{equation*}
\entire{\frac{1}{\pi}\utilde{\phi}_\nu(\lambda)+\frac{1}{2}}
\le  \mathcal{N}_{J'_\nu}(\lambda)
\le \entire{\frac{1}{\pi}\wtilde{\wtilde{\phi}}_\nu(\lambda)+\frac{1}{2}}.
\end{equation*}
\end{corollary}

\section{Proofs}\label{sec:proofs}

\subsection{Liouville's transformation and phase functions}
Consider a real-valued function $f\in C^3(a, +\infty)$,  $a\in\mathbb{R}$, with $f'(x)>0$ for $x>a$. For $t\in\mathbb{R}$, we define the functions
\begin{equation*}
\mathcal{F}_{f,t}(x):=\frac{\cos\left(f(x)-\pi t\right)}{\sqrt{f'(x)}}
\end{equation*}
and
\begin{equation}\label{eq:Qf}
\mathcal{V}_f(x):=\left(f'(x)\right)^2+\frac{1}{2}\frac{f'''(x)}{f'(x)}-\frac{3}{4}\left(\frac{f''(x)}{f'(x)}\right)^2.
\end{equation}

The following simple result is verified by a direct substitution, cf.\ \cite[Lemma 1]{EL00}. 

\begin{lemma}\label{lem:liou}
Under the above conditions, the function $\mathcal{F}_{f,t}(x)$
satisfies on $(a, +\infty)$ the Schr\"odinger equation 
\begin{equation}\label{eq:FQ}
\mathcal{F}_{f,t}''(x)+\mathcal{V}_f(x)\mathcal{F}_{f,t}(x)=0.
\end{equation}
\end{lemma}

\begin{remark} 
The relation \eqref{eq:Qf} is called  \emph{Liouville's transformation} in \cite{Hor} and \emph{Kummer's equation} in \cite[\S1]{HBRV15}.
\end{remark}

Let now $\mathcal{A}_1(x)$ and $\mathcal{A}_2(x)$ be two real linearly independent solutions of a Schr\"odinger equation 
\begin{equation}\label{eq:absSchro}
\mathcal{A}''(x)+\mathcal{P}(x)\mathcal{A}(x)=0
\end{equation}
on $(a, +\infty)$, with a real-valued potential $\mathcal{P}\in C(a, +\infty)$. We define the modulus function $\mathcal{M}(x)=\mathcal{M}_{\mathcal{A}_1,\mathcal{A}_2}(x)$ and the phase function $\Psi(x)=\Psi_{\mathcal{A}_1,\mathcal{A}_2}(x)$ for these solutions as
\begin{equation*}
\mathcal{M}(x):=\sqrt{\left(\mathcal{A}_1(x)\right)^2+\left(\mathcal{A}_2(x)\right)^2},\qquad
\mathcal{A}_1(x)+\ir \mathcal{A}_2(x)=\mathcal{M}(x)\exp\left(\ir \Psi(x)\right),
\end{equation*}
where we choose a continuous branch of $\Psi(x)=\operatorname{Arctan}\frac{\mathcal{A}_2(x)}{\mathcal{A}_1(x)}$ specified by the chosen value of 
\begin{equation*}
\lim_{x\to a^+} \Psi(x).
\end{equation*}

\begin{lemma}\label{lem:absphase} We have 
\begin{equation*}
 \Psi'(x)=\frac{\mathcal{W}\left\{\mathcal{A}_1(x), \mathcal{A}_2(x)\right\}}{\mathcal{M}(x)^2},
\end{equation*}
where $\mathcal{W}\left\{\cdot,\cdot\right\}$ is the Wronskian of the two solutions, cf.\ {\normalfont\cite[\S1]{HBRV15}}.

Assume additionally that $\Psi'(x)>0$ for $x\in (a, +\infty)$. 
Then for any $t\in\mathbb{R}$, the function
\begin{equation*}
\mathcal{A}(x):=\mathcal{F}_{\Psi, t}(x)
\end{equation*}
also satisfies \eqref{eq:absSchro}, and therefore
\begin{equation*}
\mathcal{V}_{\Psi}(x)=\mathcal{P}(x).
\end{equation*}
\end{lemma}

\begin{proof} The first statement is verified by direct differentiation: 
\begin{equation*}
\Psi'(x)=\frac{\dr}{\dr x} \arctan\frac{\mathcal{A}_2(x)}{\mathcal{A}_1(x)}=\frac{\mathcal{A}_1(x) \mathcal{A}'_2(x)-\mathcal{A}_2(x) \mathcal{A}'_1(x)}{\left(\mathcal{A}_1(x)\right)^2+ \left(\mathcal{A}_2(x)\right)^2}.
\end{equation*}
Since 
\begin{equation*}
\mathcal{F}_{\Psi, t}(x)=\frac{\cos(\pi t) \mathcal{A}_1(x)+\sin(\pi t) \mathcal{A}_2(x)}{\sqrt{\mathcal{W}}},
\end{equation*}
where $\mathcal{W}=\mathcal{W}\left\{\mathcal{A}_1(x), \mathcal{A}_2(x)\right\}$ is a constant which we assumed to be positive, the second statement also follows immediately once we have taken into account that $\mathcal{A}_1(x)$ and $\mathcal{A}_2(x)$ both satisfy \eqref{eq:absSchro}.
\end{proof}

Lemma \ref{lem:absphase} implies
\begin{lemma}
For the Bessel phase function $\theta_\nu(x)$, we have the Schr\"odinger equation \eqref{eq:FQ}, namely,
\begin{equation*}
\mathcal{F}_{\theta_\nu,t}''(x)+\mathcal{V}_{\theta_\nu}(x)\mathcal{F}_{\theta_\nu,t}(x)=0
\end{equation*}
valid in the interval $(0,+\infty)$ for all $\nu\ge 0$ and $t\in\mathbb{R}$, with
\begin{equation}\label{eq:Qtheta}
\mathcal{V}_{\theta_\nu}(x)=1-\frac{\nu^2-1/4}{x^2}.
\end{equation}
Similarly, for the Bessel derivative phase function $\phi_\nu(x)$, we also have the Schr\"odinger equation \eqref{eq:FQ}, namely,
\begin{equation*}
\mathcal{F}_{\phi_\nu,t}''(x)+\mathcal{V}_{\phi_\nu}(x)\mathcal{F}_{\phi_\nu,t}(x)=0
\end{equation*}
valid in the interval $(\nu,+\infty)$ for all $\nu\ge 0$ and $t\in\mathbb{R}$, with
\begin{equation}\label{eq:Qphi}
\mathcal{V}_{\phi_\nu}(x)=1-\frac{\nu ^2-\frac{1}{4}}{x^2}-\frac{2 \nu ^2+x^2}{\left(x^2-\nu ^2\right)^2}.
\end{equation}
\end{lemma}

\begin{proof}
We note that for all cylindrical functions $\mathcal{C}_{\nu,\tau}(x)$, and in particular for $J_\nu(x)$ and $Y_\nu(x)$, we have the Schr\"odinger equation  
\begin{equation*}
\left(\sqrt{x}\mathcal{C}_{\nu,\tau}(x)\right)''+\left(1-\frac{\nu^2-1/4}{x^2}\right)\left(\sqrt{x}\mathcal{C}_{\nu,\tau}(x)\right)=0,
\end{equation*}
see \cite[(10.13.1)]{DLMF}; it can be easily checked directly. By Lemma \ref{lem:absphase}, since $\theta_\nu(x)=\Psi_{\sqrt{x}J_\nu, \sqrt{x}Y_\nu}(x)$ has a positive derivative on $(\nu,+\infty)$, we conclude that $\mathcal{F}_{\phi_\nu,t}$ satisfies the Schr\"odinger equation with the same potential, thus proving  \eqref{eq:Qtheta}. We note that \eqref{eq:Qtheta} can be found in \cite[(10.18.16)]{DLMF}.

In the same manner, we have for $\mathcal{C}'_{\nu,\tau}(x)$, and in particular for $J'_\nu(x)$ and $Y'_\nu(x)$,
\begin{equation*}
\left(\frac{x^{3/2}}{\sqrt{x^2-\nu^2}}\,\mathcal{C}'_{\nu,\tau}(x)\right)''+\left(1-\frac{\nu ^2-\frac{1}{4}}{x^2}-\frac{2 \nu ^2+x^2}{\left(x^2-\nu ^2\right)^2}\right)\left(\frac{x^{3/2}}{\sqrt{x^2-\nu^2}}\,\mathcal{C}'_{\nu,\tau}(x)\right)=0,
\end{equation*}
which is straightforward to verify directly or deduce from \cite[(10.13.7)]{DLMF}, and the result follows from the observation that $\phi_\nu(x)=\Psi_{\frac{x^{3/2}}{\sqrt{x^2-\nu^2}}J'_\nu, \frac{x^{3/2}}{\sqrt{x^2-\nu^2}}Y'_\nu}(x)$ has a positive derivative on $(\nu,+\infty)$.
\end{proof}

\subsection{A consequence of the Sturm comparison theorem}

We state the following result which follows from the Sturm comparison theorem, and which will be used for comparing phase functions and their bounds.  This is different from a common method of comparing zeros, cf.\ \cite{Heth} and \cite{EL00}.

\begin{theorem}\label{thm:sturm} Consider, for $a\in\mathbb{R}$, two functions $g, h\in C^3(a,+\infty)$ with positive derivatives and satisfying the following conditions:
\begin{enumerate}[{\normalfont(C$_1$)}]
\item $\lim_{x\to+\infty}g(x)=\lim_{x\to+\infty}h(x)=+\infty${\normalfont;}
\item $\mathcal{V}_{g}(x) > \mathcal{V}_{h}(x)$ for all $x\in (a,+\infty)${\normalfont;}
\item there exists $b\ge a$ such that $g(x) < h(x)$ for all $x>b$.
\end{enumerate}
Then in fact $g(x) < h(x)$ for all $x>a$.
\end{theorem}

\begin{remark}\label{rem:C3prime}  
Our typical use of Theorem \ref{thm:sturm} can be illustrated by the following generic example.  

Suppose that $g(x)$ is a phase function $\Psi_{\mathcal{A}_1,\mathcal{A}_2}(x)$ of two linearly independent solutions of the Schr\"odinger equation \eqref{eq:absSchro} with a given potential $\mathcal{P}(x)$, and $h(x)$ is some conjectured upper bound for $g(x)$. Assuming that  both $g(x)$ and $h(x)$ have positive derivatives, we  can use Lemma \ref{lem:absphase} to deduce that $\mathcal{V}_{g}(x)=\mathcal{P}(x)$, and therefore verification of (C$_2$) reduces to an explicit calculation. 

In practice, we usually replace condition (C$_3$) by a stronger condition
\begin{enumerate}
\item[{\normalfont(C$'_3$)}] for some $s\in\mathbb{R}$, there exists the limit
$\lim_{x\to+\infty}x^{s}\left(h(x)-g(x)\right)$, and it is positive.
\end{enumerate}
This condition is easier to verify if the asymptotics of the phase function as $x\to+\infty$ is known. It is clear that (C$'_3$) implies (C$_3$). 
\end{remark}

\begin{proof}[Proof of Theorem \ref{thm:sturm}] Suppose that the conclusion is wrong. Then there exists $x_0>a$ for which $g(x_0) \ge h(x_0)$. Let $h_0:=h(x_0)$, $t_0:=\frac{1}{\pi}h_0-\frac{1}{2}$, and consider the following two functions,
\begin{equation*}
G(x):=\mathcal{F}_{g, t_0}(x)=\frac{\sin\left(g(x)-h_0\right)}{\sqrt{g'(x)}},\qquad 
H(x):=\mathcal{F}_{h, t_0}(x)=\frac{\sin\left(h(x)-h_0\right)}{\sqrt{h'(x)}}.
\end{equation*}
Our assumptions imply that these function are both in  $C^2(a,\infty)$, and, by Lemma \ref{lem:liou}, satisfy the differential equations
\begin{equation*}
G''(x) + \mathcal{V}_{g}(x)G(x) = 0,\qquad H''(x) + \mathcal{V}_{h}(x)H(x) = 0
\end{equation*} 
in that interval.

Due to condition (C$_1$), $H(x)$ has infinitely many zeros $\eta_k\in [x_0,\infty)$, one each time $h(\eta_k)-h_0=\pi k$, for some sequence of consecutive $k\in\mathbb Z$. Note that $\eta_0=x_0$, so that this sequence starts with $k=0$. Similarly, $G(x)$ has a zero $\gamma_k\in [x_0,\infty)$ each time $g(\gamma_k)-h_0=\pi k$ for another sequence of consecutive $k\in\mathbb Z$. Equivalently, 
\begin{equation}\label{eq:gammaeta}
\gamma_k=g^{-1}\left(\pi k+h_0\right), \qquad \eta_k=h^{-1}\left(\pi k+h_0\right).
\end{equation}

Observe that by our contradiction assumption, $g(x_0) - h_0 \ge 0$, or $\gamma_0=g^{-1}\left(h_0\right)\le x_0$, and so the first zero of $G(x)$ in $(x_0,\infty)$ is $\gamma_m$ for some $m > 0$. By condition (C$_2$) and the Sturm comparison theorem, there exists at least one zero of $G(x)$ strictly between each pair of zeros of $H(x)$ on $(x_0,\infty)$. So, $G(x)$ has a zero in $(\eta_0, \eta_1)$, and therefore we must have
\begin{equation*}
\gamma_m < \eta_1 \le \eta_m.
\end{equation*}
Similarly, there must be a zero of $G(x)$ in $\left(\eta_m, \eta_{m+1}\right)$, and therefore $\gamma_{m+1} < \eta_{m+1}$. By induction, $\gamma_k < \eta_k$ for all $k\ge m$. However, by \eqref{eq:gammaeta}
and condition (C$_3$), we must have $\gamma_k > \eta_k$ for sufficiently large $k$. This contradiction completes the proof.
\end{proof}

\subsection{Proof of {Theorem \ref{thm:thetabounds}}}\label{sec:besselproof}

We start with the upper bound. Applying Lemma \ref{lem:liou} to $f:=\wtilde{\theta}_\nu$ (which is defined and has a positive derivative on $(\nu, +\infty)$) with $t=0$, we deduce that on this interval the function $\mathcal{F}_{\wtilde{\theta}_\nu,0}$ satisfies  the Schr\"odinger equation  
\begin{equation*}
\mathcal{F}_{\wtilde{\theta}_\nu,0}''(x)+\mathcal{V}_{\wtilde{\theta}_\nu}(x)\mathcal{F}_{\wtilde{\theta}_\nu,0}(x)=0,
\end{equation*}
with the potential
\begin{equation}\label{eq:V2}
\mathcal{V}_{\wtilde{\theta}_\nu}(x)=\frac{4 x^6-12 \nu ^2 x^4+6 \nu ^2 \left(2 \nu ^2-1\right)  x^2 - \nu^4\left(4 \nu ^2-1\right)}{4 x^2 \left(x^2-\nu ^2 \right)^2}.
\end{equation}

We want to apply Theorem \ref{thm:sturm} with
\begin{equation*}
g(x)=\theta_\nu(x),\qquad
h(x)=\wtilde{\theta}_\nu(x),\qquad
a=\nu.
\end{equation*}
Condition (C$_1$) of Theorem \ref{thm:sturm} is obviously true. Verifying condition (C$_2$) is straightforward: we have, by \eqref{eq:V2} and \eqref{eq:Qtheta},
\begin{equation*}
\mathcal{V}_{\theta_\nu}(x)-\mathcal{V}_{\wtilde{\theta}_\nu}(x)=\frac{x^2+4\nu^2}{4(x^2-\nu^2)^2}>0\qquad\text{for all }x>\nu.
\end{equation*}

According to Remark \ref{rem:C3prime},  it remains to check condition  (C$'_3$). The comparison of the asymptotics \eqref{eq:thetaasympt} and \eqref{eq:overthetaasympt} yields 
\begin{equation*}
\wtilde{\theta}_\nu(x)-\theta_\nu(x)=\frac{1}{8x}+O\left(x^{-3}\right)\qquad\text{as }x\to+\infty,
\end{equation*}
and therefore (C$'_3$) holds in this case with $s=1$. Thus all conditions of Theorem \ref{thm:sturm}  are fulfilled, which proves the upper bound in \eqref{eq:thetabounds}.

Instead of proving the lower bound in \eqref{eq:thetabounds}, we first establish a weaker bound
\begin{equation}\label{eq:weakthetaboundlower}
\utilde{\theta}_\nu(x)<\theta_\nu(x).
\end{equation}
In order to prove \eqref{eq:weakthetaboundlower} we act in the same manner as above, this time applying Theorem \ref{thm:sturm} with
\begin{equation*}
g(x)=\utilde{\theta}_\nu(x),\qquad
h(x)=\theta_\nu(x),\qquad
a=\nu.
\end{equation*}
The expression for $\mathcal{V}_{\utilde{\theta}_\nu}(x)$ is given in Appendix \ref{sec:potentials}.

Condition (C$_1$) of Theorem \ref{thm:sturm} is obviously true, and verifying condition (C$_2$) requires the change of variable $x=\sqrt{\nu^2+\chi}$ with $\chi>0$, which gives, by  \eqref{eq:Qtheta} and \eqref{eq:QthetaDown1}--\eqref{eq:QthetaDown3},
\begin{gather*}
64 \chi ^5 \left(8\chi ^3+\chi ^2+6 \nu ^2 \chi +5 \nu ^4\right)^2
\left(\mathcal{V}_{\utilde{\theta}_\nu}\left(\sqrt{\nu^2+\chi}\right)-\mathcal{V}_{\theta_\nu}\left(\sqrt{\nu^2+\chi}\right)\right)=\\
1600 \chi ^9+\left(33984 \nu ^2+16\right) \chi ^8+\left(99008 \nu ^4+784 \nu ^2+1\right) \chi ^7\\
+\left(70720 \nu ^4+1696 \nu ^2+23\right) \nu ^2 \chi ^6+3 \left(1376 \nu ^2+71\right) \nu ^4 \chi ^5+\left(5200 \nu ^2+1011\right) \nu ^6 \chi ^4\\
+5 \left(400 \nu ^2+519\right) \nu ^8 \chi ^3 +3525 \nu ^{10} \chi ^2 + 2375 \nu ^{12} \chi  + 625 \nu ^{14}.
\end{gather*}
The right-hand side  and the factor in the left-hand side are obviously positive for all $\chi>0$, and thus $\mathcal{V}_{\utilde{\theta}_\nu}(x)>\mathcal{V}_{\theta_\nu}(x)$ for all $x>\nu$.

The validity of condition  (C$'_3$) follows from  comparing  \eqref{eq:thetaasympt} and \eqref{eq:underthetaasympt}, which gives
\begin{equation*}
\theta_\nu(x) - \utilde{\theta}_\nu(x) = \frac{25}{384x^3}+O\left(x^{-4}\right)\qquad\text{as }x\to+\infty,
\end{equation*}
which implies (C$'_3$) with $s=3$. This shows that Theorem \ref{thm:sturm} is applicable and proves \eqref{eq:weakthetaboundlower}. The lower bound in \eqref{eq:thetabounds} now follows as we know a priori that $\theta_\nu(x)>-\frac{\pi}{2}$ for all $x>\nu\ge 0$.

\subsection{Proof of Lemma {\ref{lem:xstar}}}

The only two non-trivial aspects are that the root $x_{\nu}^\star$ of the sextic polynomial  \eqref{eq:pnu} is greater than $\nu$, and that $z^\star_{\nu}$ is monotone decreasing in $\nu$.

To address the former, the substitution $x=\sqrt{\nu^2+\xi}$, $\xi\ge 0$, turns \eqref{eq:pnu}  into 
\begin{equation*}
p_\nu\left(\sqrt{\nu^2+\xi}\right)=8 \xi ^3-3 \xi ^2-10 \nu ^2 \xi -7 \nu ^4.
\end{equation*}
It is then elementary to check that $p_\nu(\nu)\le 0$, $p_\nu\left(\sqrt{\nu^2+\xi}\right)$ is negative near $\xi=0^+$, positive as $\xi\to+\infty$, and $\frac{\dr p_\nu\left(\sqrt{\nu^2+\xi}\right)}{\dr\xi}$ vanishes at the only positive point $\xi=\frac{1}{8}+\frac{\sqrt{80 \nu ^2+3}}{8\sqrt{3}}$, which completes the proof.

To address the latter aspect, suppose for contradiction that $\frac{\dr z^\star_\nu}{\dr\nu}$ vanishes for some $\nu=\mu>0$. 
We have
\begin{equation}\label{eq:dznudnu}
0=\left.\frac{\dr z^\star_\nu}{\dr\nu}\right|_{\nu=\mu}=\left.\frac{\dr \wtilde{\phi}_\nu\left(x^\star_\nu\right)}{\dr\nu}\right|_{\nu=\mu}
=\left.\frac{\partial \wtilde{\phi}_\nu(x)}{\partial\nu}\right|_{(\nu,x)=(\mu, x^\star_\mu)}+\wtilde{\phi}'_\mu\left(x^\star_\mu\right) \left.\frac{\dr x_\nu^\star}{\dr\nu}\right|_{\nu=\mu}.
\end{equation} 
The last term in the right-hand side (in which $'$ denotes, as usual, the derivative with respect to the argument), vanishes by the definition of $x^\star_\nu$,
and the equation reduces, after explicit evaluation of  $\frac{\partial \wtilde{\phi}_\nu(x)}{\partial\nu}$ and the substitution 
\begin{equation}\label{eq:kappa}
x^\star_\mu=\kappa \mu\qquad\text{with }\kappa>1,
\end{equation} 
to 
\begin{equation}\label{eq:mukappa}
\mu^2=\frac{23 \kappa ^2-2}{24 \left(\kappa ^2-1\right)^{5/2} \operatorname{arcsec}\kappa}.
\end{equation} 
We now recall that by \eqref{eq:pnu} we still require
\begin{equation*}
p_\mu(x^\star_\mu)=0;
\end{equation*}
substituting  \eqref{eq:kappa} and \eqref{eq:mukappa} into this gives, after some simplifications and using $\kappa>1$, the equation
\begin{equation}\label{eq:Kkappa}
\frac{23 \kappa ^4-25 \kappa ^2+2}{3 \kappa ^2 \sqrt{\kappa ^2-1} \left(3 \kappa ^2+4\right)}- \operatorname{arcsec}\kappa=0.
\end{equation}
Denote the left-hand side of \eqref{eq:Kkappa} by $K(\kappa)$. Then $\lim_{\kappa\to 1^+} K(\kappa)=0$ and
\begin{equation*}
K'(\kappa)=-\frac{16 \left(\kappa ^2-1\right)^{3/2} \left(6 \kappa ^2+1\right)}{3 \kappa ^3 \left(3 \kappa ^2+4\right)^2}
\end{equation*}
is strictly negative for $\kappa>1$. Therefore \eqref{eq:Kkappa} cannot have any solutions with $\kappa>1$, and in turn \eqref{eq:dznudnu} cannot hold for any $\mu>0$. The contradiction completes the proof.

\subsection{Proof of {Theorem \ref{thm:phibounds}}}

We act similarly to \S\ref{sec:besselproof}, starting with the lower bound in \eqref{eq:phibounds}. Since, by \eqref{eq:defunderphi}, $\utilde{\phi}_\nu(x)$ differs 
from $\wtilde{\theta}_\nu(x)$ only by a constant, Lemma \ref{lem:liou} shows that on $(\nu,+\infty)$, the function $\mathcal{F}_{\utilde{\phi}_\nu,0}$ satisfies  the Schr\"odinger equation  
\begin{equation*}
\mathcal{F}_{\utilde{\phi}_\nu,0}''(x)+\mathcal{V}_{\utilde{\phi}_\nu}(x)\mathcal{F}_{\utilde{\phi}_\nu,0}(x)=0,
\end{equation*}
with the potential 
\begin{equation}\label{eq:QphiDown}
\mathcal{V}_{\utilde{\phi}_\nu}(x)=\mathcal{V}_{\wtilde{\theta}_\nu}(x),
\end{equation}
see \eqref{eq:V2}.

We want to apply Theorem \ref{thm:sturm} with
\begin{equation*}
g(x)=\utilde{\phi}_\nu(x),\qquad
h(x)=\phi_\nu(x),\qquad
a=\nu.
\end{equation*}
By \eqref{eq:Qphi}, \eqref{eq:QphiDown}, and \eqref{eq:V2},
\begin{equation*}
\mathcal{V}_{\utilde{\phi}_\nu}(x)-\mathcal{V}_{\phi_\nu}(x)=\frac{4 \nu ^2+3 x^2}{4 \left(x^2-\nu ^2\right)^2}>0\qquad\text{for all }x>\nu,
\end{equation*}
thus verifying condition (C$_2$) of Theorem \ref{thm:sturm}. Further, by \eqref{eq:phiasympt} and \eqref{eq:underphiasympt},
\begin{equation*}
\phi_\nu(x)-\utilde{\phi}_\nu(x)=\frac{3}{8x}+O\left(x^{-2}\right)\qquad\text{as }x\to+\infty,
\end{equation*}
confirming that condition (C$'_3$) holds with $s=1$. Thus,  Theorem \ref{thm:sturm} is applicable and yields the lower bound in \eqref{eq:phibounds}. 

Before proving the upper bound in \eqref{eq:phibounds} in full generality, we will prove it for sufficiently large $x$ by showing that
\begin{equation}\label{eq:phiboundxlarge}
\phi_\nu(x)<\wtilde{\phi}_\nu(x)\qquad\text{for }x>x_\nu^\star. 
\end{equation}
Once this is done, the full bound follows easily: firstly, since $\phi'_\nu\left(x_\nu^\star\right)>0=\wtilde{\phi}'_\nu\left(x_\nu^\star\right)$, the strict inequality also holds for $x=x_\nu^\star$; the rest follows since $\phi_\nu(x)<\phi_\nu\left(x_\nu^\star\right)$ and  $\wtilde{\phi}_\nu(x)>\wtilde{\phi}_\nu\left(x_\nu^\star\right)=z_\nu^\star$ for $x\in\left(\nu,x_\nu^\star\right)$.

We will prove \eqref{eq:phiboundxlarge} by applying Theorem \ref{thm:sturm} with 
\begin{equation*}
g(x)=\phi_\nu(x),\qquad
h(x)=\wtilde{\phi}_\nu(x),\qquad
a=x_\nu^\star.
\end{equation*}
Since $\wtilde{\phi}'_\nu(x)>0$ for $x\in\left(x_\nu^\star,+\infty\right)$ by Lemma \ref{lem:xstar}, we can apply Lemma \ref{lem:liou} on this interval to deduce that the function $\mathcal{F}_{\wtilde{\phi}_\nu,0}$ satisfies there the Schr\"odinger equation  
\begin{equation*}
\mathcal{F}_{\wtilde{\phi}_\nu,0}''(x)+\mathcal{V}_{\wtilde{\phi}_\nu}(x)\mathcal{F}_{\wtilde{\phi}_\nu,0}(x)=0,
\end{equation*}
where the potential  $\mathcal{V}_{\wtilde{\phi}_\nu}(x)$ is given by \eqref{eq:QphiUp}--\eqref{eq:QphiUp4}.

In order to verify condition (C$_2$), we set
\begin{equation*}
\delta_\nu(x):=4096 x^2 \left(x^2-\nu ^2\right)^{10} \left(\wtilde{\phi}'_\nu(x)\right)^2 \left(\mathcal{V}_{\phi_\nu}(x)-\mathcal{V}_{\wtilde{\phi}_\nu}(x)\right).
\end{equation*}
Since the first three factors in the definition of $\delta_\nu(x)$ are positive for $x>x_\nu^\star$, it is enough to show that 
\begin{equation}\label{eq:deltaineq}
\delta_\nu(x)>0\qquad\text{for all }x>x_\nu^\star
\end{equation}
to ensure that (C$_2$) holds. Explicit computation gives
\begin{equation*}
\begin{split}
\delta_\nu(x)=&\hphantom{+\,}4032 x^{18}+16 \left(1208 \nu ^2+27\right) x^{16}-\left(158656 \nu ^4+2640 \nu ^2+81\right) x^{14}\\
&+16 \nu ^2 \left(20272 \nu ^4+701 \nu ^2-27\right) x^{12}-16 \nu ^4 \left(12796 \nu ^4+1415 \nu ^2+54\right) x^{10}\\
& -64 \nu ^6 \left(2338 \nu ^4-337 \nu ^2+12\right) x^8+64 \nu ^8 \left(4543 \nu ^4-164 \nu ^2-4\right) x^6\\
&-3584 \nu ^{12} \left(41 \nu ^2-1\right) x^4+1024 \nu ^{14} \left(17 \nu ^2-1\right) x^2+4096 \nu ^{18},
\end{split}
\end{equation*}
and we will verify \eqref{eq:deltaineq} separately in two cases.

\begin{description}
\item[Case $\nu=0$. ] Since $x_0^\star=\sqrt{\frac{3}{8}}$ and 
\begin{equation*}
\left.x^{-14} \delta_0(x)\right|_{x=\sqrt{\frac{3}{8}+\xi}}=72 \left(56 \xi ^2+48 \xi +9\right)>0\qquad\text{for }\xi>0,
\end{equation*}  
\eqref{eq:deltaineq} follows.
\item[Case $\nu>0$. ] We introduce the new variable $\zeta>0$ via $x=\nu(1+\zeta)$, and note that
\begin{equation*}
\begin{split}
\nu^{-14}\delta_\nu(x)&=64 \zeta ^6 (\zeta +2)^6 \left(63 \zeta ^6+378 \zeta ^5+1625 \zeta ^4+3980 \zeta ^3+5681 \zeta ^2+4410 \zeta +1463\right) \nu ^4\\
&+16 \zeta ^3 (\zeta +1)^2 (\zeta +2)^3 \\
&\ \times\left(27 \zeta ^8+216 \zeta ^7+672 \zeta ^6+1008 \zeta ^5+998 \zeta ^4+1304 \zeta ^3+1672 \zeta ^2+1120 \zeta +343\right) \nu ^2\\
&-(\zeta +1)^6 \left(3 \zeta ^2+6 \zeta +7\right)^4,
\end{split}
\end{equation*}
and, since the coefficient of the $\nu^2$-term above is always positive, 
\begin{equation}\label{eq:deltaineq1}
\begin{split}
\nu^{-14}\delta_\nu(x)&\ge 64 \zeta ^6 (\zeta +2)^6 \left(63 \zeta ^6+378 \zeta ^5+1625 \zeta ^4+3980 \zeta ^3+5681 \zeta ^2+4410 \zeta +1463\right) \nu ^4
\\&-(\zeta +1)^6 \left(3 \zeta ^2+6 \zeta +7\right)^4.
\end{split}
\end{equation}
On the other hand, we have 
\begin{equation*}
x>x_\nu^\star \iff \nu^{-4} p_\nu(\nu (1+\zeta))>0\\
\iff \nu ^4\geq \frac{(\zeta +1)^4 (3 \zeta  (\zeta +2)+7)^2}{64 \zeta ^6 (\zeta +2)^6}.
\end{equation*}
Taking this into account in \eqref{eq:deltaineq1}, we deduce that for $x>x_\nu^\star$,
\begin{equation*}
\nu^{-14}\delta_\nu(x)\ge 2 (\zeta +1)^4 (3 \zeta  (\zeta +2)+7)^2 (\zeta  (\zeta +2) (\zeta  (\zeta +2) (27 \zeta  (\zeta +2)+409)+1057)+707) >0,
\end{equation*}
proving \eqref{eq:deltaineq} in this case as well.
\end{description}

To finish the proof of \eqref{eq:phiboundxlarge} using Theorem \ref{thm:sturm}, it remains to check condition (C$'_3$), for which we use  \eqref{eq:phiasympt} and \eqref{eq:overphiasympt} to deduce
\begin{equation*}
\wtilde{\phi}_\nu(x)-\phi_\nu(x)=\frac{21}{128x^3}+O\left(x^{-4}\right)\qquad\text{as }x\to+\infty,
\end{equation*}
confirming that condition (C$'_3$) holds with $s=3$.

\section{Derivatives of ultraspherical Bessel functions}\label{sec:ultra}

\subsection{Setup III}
\begin{figure}[ht]
\centering
\includegraphics{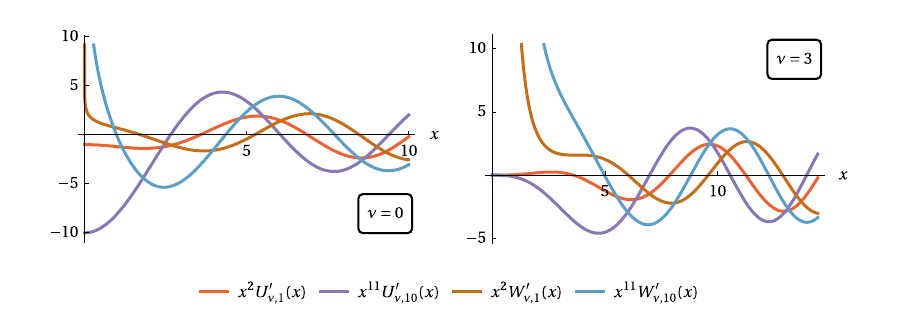}
\caption{Plots of $x^{\eta+1}U'_{\nu, \eta}(x)$ and $x^{\eta+1}W'_{\nu, \eta}(x)$.\label{fig:UWprime}}
\end{figure}

Let $\nu\ge 0$, $\eta\in\mathbb{R}$, and let
\begin{equation*}
U_{\nu, \eta}(x):=x^{-\eta}J_{\nu}(x),\qquad W_{\nu, \eta}(x):=x^{-\eta}Y_{\nu}(x)
\end{equation*}
denote the \emph{ultraspherical Bessel functions}. We will be interested in their derivatives
\begin{equation*}
U'_{\nu, \eta}(x)=x^{-\eta-1}\left(x J'_{\nu}(x) - \eta J_\nu(x)\right),\qquad W'_{\nu, \eta}(x)=x^{-\eta-1}\left(x Y'_{\nu}(x) - \eta Y_\nu(x)\right),
\end{equation*}
in practise usually omitting the factor $x^{-\eta-1}$,  see Figure \ref{fig:UWprime}, and denote by 
\begin{equation*}
u'_{\nu, \eta, k}\qquad\text{and}\qquad w'_{\nu, \eta, k}
\end{equation*}
the $k$th positive zero of $U'_{\nu, \eta}(x)$ and $W'_{\nu, \eta}(x)$, respectively, with the exception in case $\eta=\nu$, 
\begin{equation*}
u'_{\nu, \nu, 1}:=0.
\end{equation*}
Of course, for $\eta=0$ we return to the usual derivatives of Bessel functions and their zeros.

The function $U'_{\nu, \eta}(x)$ and, specifically, its zeros $u'_{\nu, \eta, k}$, play an important role in spectral geometry. In particular, for $d\ge 3$, the numbers
\begin{equation*}
\left(u'_{m+\frac{d}{2}-1, \frac{d}{2}-1, k}\right)^2,\qquad m\in\mathbb{N}\cup\{0\}, 
\end{equation*}
are the eigenvalues of multiplicity
\begin{equation*}
\binom{m+d-1}{d-1}-\binom{m+d-3}{d-1}
\end{equation*}
of the Neumann Laplacian on the unit ball in $\mathbb{R}^d$, see \cite{AB93} and \cite[\S1.2.3]{LMP}. Also,
\begin{equation*}
\left(u'_{m, \eta, k}\right)^2,\qquad m\in\mathbb{N}\cup\{0\},
\end{equation*}
are the non-negative eigenvalues of the Robin Laplacian in the unit disk in $\mathbb{R}^2$ with the Robin parameter $-\eta$, which are taken with multiplicity one if $m=0$, and multiplicity two otherwise, see \cite{BFK}.

For further properties of the derivatives of the ultraspherical Bessel functions and of their roots, including some bounds, see, e.g., \cite{IS88}, \cite{LS}, and \cite{L99}. In turn, we are primarily interested in the corresponding phase function which we define and discuss in the next section.

\subsection{Phase function of ultraspherical Bessel derivatives}
As for the standard Bessel functions and their derivatives, we introduce the \emph{modulus} and \emph{phase functions} of the ultrasperical Bessel derivatives, $L_{\nu,\eta}, \psi_{\nu,\eta}: (0,+\infty)\to\mathbb{R}$ defined by 
\begin{equation*}
L_{\nu,\eta}(x)=\sqrt{\left(U'_{\nu,\eta}(x)\right)^2+\left(W'_{\nu,\eta}(x)\right)^2},\qquad U'_{\nu,\eta}(x)+\ir W'_{\nu,\eta}(x)=L_{\nu,\eta}(x)\left(\cos\psi_{\nu,\eta}(x)+\ir \sin\psi_{\nu,\eta}(x)\right),
\end{equation*}
where we choose a continuous branch of $\psi_{\nu,\eta}(x)$ with the initial condition 
\begin{equation*}
\lim_{x\to 0^+}\psi_{\nu,\eta}(x)=\frac{\pi}{2}.
\end{equation*}

Set, for $\nu\ge 0$, \footnote{This is a convenient notation as we will not need to consider $\mu_{\nu,\eta}$ when $\mu_{\nu,\eta}^2$ is negative.}
\begin{equation*}
\mu^2_{\nu,\eta}=\mu^2:=\nu^2-\eta^2,
\end{equation*}
and
\begin{equation*}
I_{\nu,\eta} :=\begin{cases} 
\left(\mu_{\nu, \eta},+\infty\right)&\qquad\text{if }|\eta|<\nu,\\
\left(0,+\infty\right)&\qquad\text{if }|\eta|\ge \nu.
\end{cases}
\end{equation*}
We have

\begin{lemma}\label{lem:psimonotone} 
Let $\nu\ge 0$, $\eta\in\mathbb{R}$, $x>0$. 
If $\mu^2_{\nu,\eta}>0$, then the function $\psi'_{\nu, \eta}(x)$ has a single zero at $x=\mu_{\nu,\eta}$ and is positive for $x>\mu_{\nu,\eta}$, otherwise its derivative is positive on the positive real line. Moreover, on the interval $I_{\nu,\eta}$, the function $\mathcal{F}_{\psi_{\nu, \eta},0}(x)$ satisfies the Schr\"odinger equation \eqref{eq:FQ} with the potential
\begin{equation*}
\mathcal{V}_{\psi_{\nu, \eta}}(x)=1 - \frac{\nu ^2-\frac{1}{4}}{x^2} + \frac{2 (1-\eta )}{x^2-\mu_{\nu,\eta}^2} - \frac{3 x^2}{\left(x^2-\mu_{\nu,\eta}^2\right)^2}.
\end{equation*}
\end{lemma}

\begin{proof} 
A direct differentiation shows that 
\begin{equation}\label{eq:psiprimeeq}
\psi'_{\nu, \eta}(x)=\left(\arctan\frac{W'_{\nu,\eta}(x)}{U'_{\nu,\eta}(x)}\right)'=\frac{2\left(x^2-\mu^2_{\nu,\eta}\right)}{\pi x^{2\eta+3}L_{\nu,\eta}^2(x)},
\end{equation}
and the first statement then follows immediately since $L_{\nu,\eta}(x)$ does not vanish.

Set, for $x\in I_{\nu,\eta}$,
\begin{equation*}
\mathcal{A}_1(x):=\frac{x^{\eta+3/2}}{\sqrt{x^2-\mu^2_{\nu,\eta}}}U'_{\nu, \eta}(x),\qquad
\mathcal{A}_2(x):=\frac{x^{\eta+3/2}}{\sqrt{x^2-\mu^2_{\nu,\eta}}}W'_{\nu, \eta}(x).
\end{equation*}
A straightforward check shows that we have
\begin{equation*}
\mathcal{A}''_j(x)+\left(1 - \frac{\nu ^2-\frac{1}{4}}{x^2} + \frac{2 (1-\eta )}{x^2-\mu_{\nu,\eta}^2} - \frac{3 x^2}{\left(x^2-\mu_{\nu,\eta}^2\right)^2}\right)\mathcal{A}_j(x)=0,\qquad j=1,2,
\end{equation*}
and the second statement is obtained by Lemma \ref{lem:liou}.
\end{proof}

Lemma \ref{lem:psimonotone}  ensures that the inverse function $\left(\psi_{\nu, \eta}\right)^{-1}$ is well-defined on $\left[\psi_{\nu, \eta}\left(\mu_{\nu,\eta}\right),+\infty\right)$ if $|\eta|<\nu$, or on $[0,+\infty)$ otherwise, and 
\begin{equation}\label{eq:uwprimek}
u'_{\nu, \eta,k}=\left(\psi_{\nu, \eta}\right)^{-1}\left(\pi\left(k-\frac{1}{2}\right)\right),\qquad
w'_{\nu, \eta,k}=\left(\psi_{\nu, \eta}\right)^{-1}(\pi k).
\end{equation}

We can establish the asymptotics of $\psi_{\nu, \eta}(x)$ for large $x$.

\begin{lemma} We have
\begin{equation}\label{eq:psiasympt}
\begin{split}
\psi_{\nu, \eta}(x)&=x - \frac{\pi}{4}(2 \nu -1) +\frac{4 \nu ^2+3 + 8 \eta}{8 x}\\
&+\frac{16 \nu ^4 +(192 \eta +184) \nu ^2  - 63 -128 \eta ^3-192 \eta ^2-144 \eta}{384 x^3}+O\left(x^{-5}\right)\\
&\qquad\qquad\qquad\qquad\qquad\qquad\qquad\qquad\qquad\qquad\text{as }x\to+\infty.
\end{split}
\end{equation}
\end{lemma}

\begin{proof} 
In essence, we emulate the approach of \cite[\S 2]{HBRV15} to the asymptotics of the standard Bessel phase function. 
Using standard asymptotics of the Bessel functions and their derivatives \cite[\S 10.17]{DLMF}, we get first
\begin{equation*}
x^{2\eta} L_{\nu,\eta}^2(x)
=\frac{2}{\pi  x}
+\frac{8 \eta ^2+8 \eta -4 \nu ^2+3}{4 \pi  x^3}
+\frac{\left(4 \nu ^2-1\right) \left(16 \eta ^2+48 \eta -4 \nu ^2+45\right)}{64 \pi  x^5}
+O\left(x^{-7}\right)
\end{equation*}
as $x\to+\infty$. Computing the reciprocal of this series, substituting it into the right-hand side of the differential equation \eqref{eq:psiprimeeq}, and integrating term by term gives \eqref{eq:psiasympt} except for the $O(1)$-term which, at this stage, is an unknown constant of integration $C$. Then, in the leading terms,
\begin{equation*}
U'_{\nu,\eta}(x)=L_{\nu,\eta}(x)\cos \psi_{\nu, \eta}(x) \sim \sqrt{\frac{2}{\pi }} x^{-\eta -\frac{1}{2}} \cos \left(x+C\right),
\end{equation*} 
which should match the standard asymptotics
\begin{equation*}
U'_{\nu,\eta}(x) \sim \sqrt{\frac{2}{\pi }} x^{-\eta -\frac{1}{2}} \cos \left(x-\frac{\pi}{4} (2 \nu -1)\right),
\end{equation*}
giving
\begin{equation*}
C=-\frac{\pi}{4} (2 \nu -1).
\end{equation*}
\end{proof}

We are going to demonstrate that the uniform bounds on $\psi_{\nu,\eta}(x)$ can be obtained in the same manner as we have done for $\theta_\nu(x)$ and $\phi_\nu(x)$, subject to some restrictions on $x$. Generally speaking, such bounds depend upon the signs of $\eta$ and of $\mu^2_{\nu, \eta}$, so for brevity we from now on restrict ourselves to the case
\begin{equation}\label{eq:nueta}
\nu>\eta>0
\end{equation}
(which implies $\mu^2_{\nu, \eta} >0$) and to the lower bound on $\psi_{\nu,\eta}(x)$. See also Remark \ref{lem:mu0} for the case $\nu=\eta>0$.

\subsection{Definitions and properties of the auxiliary functions III}\label{subseq:auxIII}

Assuming \eqref{eq:nueta}, we from now on use the shorthand $\mu$ for $\mu_{\nu,\eta}=\sqrt{\nu^2-\eta^2}$. 
We will keep using indices $\{\nu,\eta\}$ for some quantities but $\{\mu,\eta\}$ for some others where it simplifies the presentation, keeping the relation above in mind.
 
We define, for $x>\mu$, the function
\begin{equation}\label{eq:psidown} 
\utilde{\psi}_{\nu,\eta}(x):=
\sqrt{x^2-\mu ^2}
-\left(\frac{\eta ^2}{2 \mu }+\mu \right) \arccos\frac{\mu }{x} 
+\frac{\eta }{\sqrt{x^2-\mu ^2}} 
+ \frac{\pi }{4}  \left(\frac{\eta ^2}{\mu }+2 (\mu -\nu )+1\right).
\end{equation}

\begin{remark} Similarly to Remarks \ref{rem:asDebye} and \ref{rem:asDebyeprime}, the function $\utilde{\psi}_{\nu, \eta}(x)$ coincides with the two-term asymptotic expansions of $\psi_{\nu, \eta}(x)$ in the Debye regime (R$_\text{iii}$) deduced using the methods of \cite{Olv74}.
\end{remark}

It is obvious that
\begin{equation*}
\lim_{x\to\mu^+}\utilde{\psi}_{\nu,\eta}(x)=\lim_{x\to+\infty}\utilde{\psi}_{\nu,\eta}(x)=+\infty,
\end{equation*}
and easily checked that
\begin{equation}\label{eq:psiDownasympt}
\utilde{\psi}_{\nu, \eta}(x) =x - \frac{\pi}{4}(2 \nu -1) +\frac{\nu ^2+2 \eta}{2 x}+O\left(x^{-3}\right)\qquad\text{as }x\to+\infty.
\end{equation} 
Further, it is easy to check that its derivative
\begin{equation*}
\utilde{\psi}'_{\nu,\eta}(x)=\frac{2 x^4-x^2 \left(\eta  (\eta +2)+4 \mu ^2\right)+\mu ^2 \left(\eta ^2+2 \mu ^2\right)}{2 x(x^2-\mu^2 )^{3/2}}
\end{equation*} 
vanishes at the point
\begin{equation*}
x_{\mu,\eta}^\#:=\frac{1}{2}\sqrt{4 \mu ^2 +\eta(\eta +2) + \sqrt{\eta^2(\eta +2)^2+16 \mu ^2\eta}}\ >\mu,
\end{equation*}
and is positive for $x>x_{\mu,\eta}^\#$. 

Therefore, the inverse function $\left(\utilde{\psi}_{\nu,\eta}\right)^{-1}:\left(\utilde{\psi}_{\nu,\eta}\left(x_{\mu,\eta}^\#\right),+\infty\right)\to \left(x_{\mu,\eta}, +\infty\right)$ is well-defined.

Set additionally, for $\mu\ge 0$, $\eta>0$,
\begin{equation}\label{eq:r} 
\begin{split}
r_{\mu,\eta}(x):=
&\hphantom{+\,}
12 x^{14}\\
&+x^{12} \left(-44 \mu ^2 + 4 \eta  \left(\eta ^3+4 \eta ^2-5 \eta -18\right)\right)\\
&+x^{10} \left(40 \mu ^4+8 \left(-3 \eta ^3-10 \eta ^2+6 \eta +5\right) \eta  \mu ^2-\eta ^2 (\eta +2)^2 \left(4 \eta ^2+8 \eta -3\right)\right)\\
&+x^8 \left(40 \mu ^6+4 \left(15 \eta ^3+40 \eta ^2+2 \eta +70\right) \eta  \mu ^4\right.\\
&\qquad\left. +\left(20 \eta ^4+96 \eta ^3+147 \eta ^2+96 \eta +52\right) \eta ^2 \mu ^2+(\eta +2)^4 \eta ^4\right)\\
&+x^6 \left(-100 \mu ^8-8 \eta  \left(10 \eta ^3+20 \eta ^2+14 \eta +45\right) \mu ^6\right.\\
&\qquad\left. -\eta ^2 \left(40 \eta ^4+144 \eta ^3+171 \eta ^2+116 \eta +80\right) \mu ^4-4 \eta ^5 (\eta +2)^3 \mu ^2\right)\\
&+x^4 \left(68 \mu ^{10} +4 \left(15 \eta ^3+20 \eta ^2+27 \eta +20\right) \eta  \mu ^8\right.\\
&\qquad\left.+\left(40 \eta ^4+96 \eta ^3+81 \eta ^2+24 \eta +16\right) \eta ^2 \mu ^6 + 6 (\eta +2)^2 \eta ^6 \mu ^4\right)\\
&+x^2 \left(-16 \mu ^{12}-8 \eta  \left(3 \eta ^3+2 \eta ^2+4 \eta -4\right) \mu ^{10}\right.\\
&\qquad\left.-4 \eta ^3 \left(5 \eta ^3+6 \eta ^2+3 \eta -4\right) \mu ^8-4 \eta ^7 (\eta +2) \mu ^6\right)\\
&+\left(4 \eta ^4 \mu ^{12} +4 \eta ^6 \mu ^{10} + \eta ^8 \mu ^8\right),
\end{split}
\end{equation}
and let $x^@_{\mu,\eta}$ be the greatest positive real root of \eqref{eq:r}. It exists due to 

\begin{lemma} For any  $\mu\ge 0$, $\eta>0$,
\begin{equation*}
x^@_{\mu,\eta} > x^\#_{\mu,\eta}.
\end{equation*}
\end{lemma}

\begin{proof} Since the polynomial \eqref{eq:r} has a positive leading term $12x^{14}$, it is sufficient to show that $r_{\mu,\eta}\left(x^\#_{\mu,\eta}\right)<0$. 
With the shorthand $\rho:=\sqrt{\eta^2 (\eta +2)^2 + 16 \mu^2\eta}>0$, we get
\begin{equation*}
\begin{split}
r_{\mu,\eta}\left(x^\#_{\mu,\eta}\right)&=-\frac{3}{16} \eta ^2 \left(\eta  (\eta +2)^2+16 \mu ^2\right)\\
&\times \left(\eta ^9+10 \eta ^8+\eta ^7 \left(4 \mu ^2+\rho +40\right)+4 \eta ^6 \left(13 \mu ^2+2 \rho +20\right)\right.\\
&\qquad\left.+4 \eta ^5 \left(\mu ^4+\mu ^2 (\rho +54)+6 \rho +20\right)+4 \eta ^4 \left(22 \mu ^4+\mu ^2 (9 \rho +92)+8 (\rho +1)\right)\right.\\
&\qquad\left.+4 \eta ^3 \left(\mu ^4 (\rho +96)+8 \mu ^2 (3 \rho +7)+4 \rho \right)+16 \eta ^2 \left(3 \mu ^6+\mu ^4 (3 \rho +28)+5 \mu ^2 \rho \right)\right.\\
&\qquad\left.+32 \eta  \left(7 \mu ^6+3 \mu ^4 \rho \right)+16 \mu ^6 \rho \right)<0 
\end{split}
\end{equation*}
as required.
\end{proof}

\subsection{Main results III: bounding the phase and zeros of derivatives of  ultraspherical Bessel functions}

\begin{theorem}\label{thm:psibound}  For every $\nu>\eta>0$ and every $x>x_{\mu,\eta}^@$, with $\mu=\mu_{\nu,\eta}$, we have
\begin{equation*}
\utilde{\psi}_{\nu,\eta}(x)<\psi_{\nu,\eta}(x).
\end{equation*}
\end{theorem}

Recalling \eqref{eq:uwprimek}, we deduce

\begin{corollary}\label{cor:uwprimebounds} 
For $\nu>\eta>0$, we have the following bounds, valid for all $k\in\mathbb{N}$ such that the argument of the inverse function exceeds $\utilde{\psi}_{\nu,\eta}\left(x_{\mu,\eta}^@\right)$.
\begin{enumerate}[{\normalfont(i)}]
\item 
\begin{equation*}
u'_{\nu,\eta, k}<\left(\utilde{\psi}_\nu\right)^{-1}\left(\pi\left(k-\frac{1}{2}\right)\right).
\end{equation*}
\item
\begin{equation*}
w'_{\nu,\eta, k}<\left(\utilde{\psi}_\nu\right)^{-1}\left(\pi k\right).
\end{equation*}
\end{enumerate}
\end{corollary}

\begin{remark}\label{lem:mu0} Formally, Theorem \ref{thm:psibound} and Corollary \ref{cor:uwprimebounds} are inapplicable when $\eta=\nu$, and therefore $\mu=0$, since in that case the definition \eqref{eq:psidown}  does not make sense. If, however, we set, for $\nu>0$ and $x>0$,
\begin{equation*}
\utilde{\psi}_{\nu,\nu}(x):=\lim_{\eta\to\nu^-}\utilde{\psi}_{\nu,\eta}(x)=
x - \frac{\pi}{4} (2\nu -1 )+\frac{\nu  (\nu +2)}{2 x},
\end{equation*}
both results become valid for $\eta=\nu$ without further modifications.
\end{remark}

\begin{figure}[ht]
\centering
\includegraphics{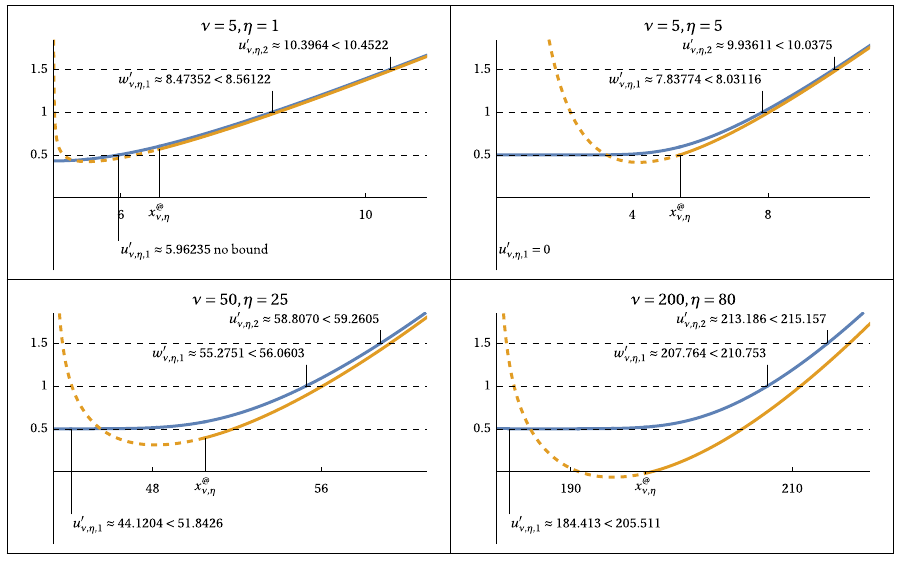}
\caption{All plots show the scaled phase function $\frac{1}{\pi}\psi_\nu(x)$ (blue), and its scaled bound $\frac{1}{\pi}\utilde{\psi}_\nu(x)$ (orange). The quality of the bounds for $u'_{\nu,\eta,k}$ and  $w'_{\nu,\eta,k}$ increases with $k$.\label{fig:comparisonU}}
\end{figure}

For an illustration of  the results in Theorem  \ref{thm:psibound}, Corollary \ref{cor:uwprimebounds}, and Remark \ref{lem:mu0}, see Figure \ref{fig:comparisonU}.

\begin{proof}[Proof of Theorem \ref{thm:psibound}]
We will prove \eqref{eq:phiboundxlarge} by applying Theorem \ref{thm:sturm} with 
\begin{equation*}
g(x)=\utilde{\psi}_{\nu,\eta}(x),\qquad
h(x)=\psi_{\nu,\eta}(x),\qquad
a=x^@_{\mu,\eta}.
\end{equation*}
Both functions $g$ and $h$ have positive derivatives on $(a, +\infty)$ by Lemma \ref{lem:psimonotone} and analysis at the beginning of \S\ref{subseq:auxIII} (since $x^@_{\mu,\eta}>x^\#_{\mu,\eta}>\mu$) and tend to infinity at infinity.
We can apply Lemma \ref{lem:liou} on this interval to deduce that the function $\mathcal{F}_{\utilde{\psi}_{\nu,\eta},0}$ satisfies there the Schr\"odinger equation  
\begin{equation*}
\mathcal{F}_{\utilde{\psi}_{\nu,\eta},0}''(x)+\mathcal{V}_{\utilde{\psi}_{\nu,\eta}}(x)\mathcal{F}_{\utilde{\psi}_{\nu,\eta},0}(x)=0,
\end{equation*}
with the potential  $\mathcal{V}_{\utilde{\psi}_{\nu,\eta}}(x)$ given by \eqref{eq:QpsiDown}--\eqref{eq:QpsiDown6}.

An explicit calculation yields 
\begin{equation*}
16 x^4 (x^2-\mu^2)^6 \left(\utilde{\psi}'_{\nu,\eta}(x)\right)^2 \left(\mathcal{V}_{\utilde{\psi}_{\nu,\eta}}(x)-\mathcal{V}_{\psi_{\nu,\eta}}(x)\right)
=r_{\mu,\eta}(x),
\end{equation*}
where the right-hand side is given by \eqref{eq:r}. Since $r_{\mu,\eta}(x)>0$ for $x>a$, we have verified condition (C$_2$) of Theorem \ref{thm:sturm}. 

Condition (C$'_3$) holds with $s=1$ since by \eqref{eq:psiasympt} and \eqref{eq:psiDownasympt},
\begin{equation*}
\psi_{\nu, \eta}(x)-\utilde{\psi}_{\nu, \eta}(x)=\frac{3}{8 x}+O\left(x^{-3}\right)\qquad\text{as }x\to+\infty.
\end{equation*} 
\end{proof}

\section{Benchmarking and conclusions}\label{sec:benchmark}

In oder to set up comparison with previous bounds, we denote
\begin{equation*}
\begin{split}
A_\nu^{(1)}(\beta)&:=\beta,\\
A_\nu^{(2)}(\beta)&:=\beta-\frac{4\nu^2-1}{8\beta},\\
A_\nu^{(3)}(\beta)&:=\beta-\frac{4\nu^2-1}{8\beta}-\frac{4(4\nu^2-1)(28\nu^2-31)}{3(8\beta)^3}.
\end{split}
\end{equation*}
and let
\begin{equation*}
\beta_{\nu, k}:= \pi\left(k+\frac{\nu}{2}-\frac{1}{4}\right).
\end{equation*}
The $i$-term McMahon asymptotic expansion for the large zeros of the Bessel function of a \emph{fixed} order $\nu\ge 0$, is given, for  $i=1,2,3$, by the formula 
\begin{equation}\label{eq:jnukasympt}
j_{\nu,k}=A_\nu^{(i)}\left(\beta_{\nu,k}\right)+O\left(\beta_{\nu, k}^{-2i+1}\right)\qquad\text{as }k\to+\infty,
\end{equation}
see \cite{McM}, \cite[\S15.53]{Wat},  \cite[(10.21.19)]{DLMF}, and also \cite[Appendix A]{Nem21} for further terms. The remainder term may be alternatively written as $O\left(k^{-2(i-1)}\right)$.
The more general asymptotic formulae for $c_{\nu,\tau,k}$ are obtained in the same manner by replacing $\beta_{\nu,k}$ with $\pi\left(k+\frac{\nu}{2}+\tau-\frac{5}{4}\right)$.

We have Hethcote's upper bound \cite{Heth}: 
\begin{equation*}
j_{\nu, k}\le
\begin{cases}
A_\nu^{(2)}(\beta_{\nu,k})\quad&\text{if }\nu\le \frac{1}{2},\\
A_\nu^{(1)}(\beta_{\nu,k})\quad&\text{if }\nu> \frac{1}{2},
\end{cases}
\quad=:\wtilde{h}_{\nu,k}\qquad\text{for  all }\nu\ge 0\text{ and }k\in\mathbb{N}.
\end{equation*}
In the same paper, Hethcote also proves the lower bound (but only valid for small $\nu$),
\begin{equation*}
j_{\nu, k}\ge H_\nu^{(1)}(\beta_{\nu,k})=:{\utilde{h}}_{\nu,k}\qquad\text{for  all }0\le \nu\le \frac{1}{2}\text{ and }k\in\mathbb{N}.
\end{equation*}

Some improvements of Hethcote's bounds were obtained in \cite{EL00}.  Namely, they proved the upper bound
\begin{equation*}
j_{\nu, k}\le
\begin{cases}
A_\nu^{(2)}(\beta_{\nu,k})\quad&\text{if }\nu\le \frac{1}{2},\\
A_\nu^{(3)}(\beta_{\nu,k})\quad&\text{if }\nu> \frac{1}{2},
\end{cases}
\quad=:\wtilde{\ell}_{\nu,k}\qquad\text{for  all }\nu\ge 0\text{ and }k\in\mathbb{N}.
\end{equation*}
and the lower bound (for a slightly wider range of $\nu$s)
\begin{equation*}
j_{\nu, k}\ge
\begin{cases}
A_\nu^{(3)}(\beta_{\nu,k})\quad&\text{if }\nu\le \frac{1}{2},\\
A_\nu^{(2)}(\beta_{\nu,k})\quad&\text{if }\frac{1}{2}<\nu<\sqrt{\frac{31}{28}},
\end{cases}
\quad=:\utilde{\ell}_{\nu,k}\qquad\text{for  all }0\le \nu< \sqrt{\frac{31}{28}}\text{ and }k\in\mathbb{N}.
\end{equation*}
For further improvements of these bounds when $\nu\le \frac{1}{2}$, see \cite{Nem21}.

We note that we always have (when the bounds exist)
\begin{equation}\label{eq:jnukHEL}
\utilde{h}_{\nu,k}\le\utilde{\ell}_{\nu,k}<j_{\nu,k}<\wtilde{\ell}_{\nu,k}\le\wtilde{h}_{\nu,k}.
\end{equation}

Alternative bounds, based on the asymptotics of $j_{\nu,k}$ in the regime (R$_\text{ii}$), were given by Qu and Wong \cite{QW}. Let
\begin{equation*}
(0>)a_1>a_2>\dots
\end{equation*}
denote the sequence of negative zeros of the Airy function $\operatorname{Ai}(x)$. Then
\begin{equation}\label{eq:jnukpQW}
\begin{split}
\utilde{q}_{\nu,k}&:=\nu-a_k\left(\frac{\nu}{2}\right)^{1/3}<j_{\nu,k}<\nu-a_k\left(\frac{\nu}{2}\right)^{1/3}+\frac{3}{20}a_k^2\left(\frac{\nu}{2}\right)^{-1/3}=:\wtilde{q}_{\nu,k}\\
&\qquad\qquad\qquad\qquad\qquad\qquad\qquad\qquad\qquad\qquad\qquad\qquad\text{for  all }\nu>0 \text{ and }k\in\mathbb{N}.
\end{split}
\end{equation}

Switching to existing bounds for $j'_{\nu,k}$, we note that although the McMahon asymptotic expansion for these zeros for a fixed $\nu$ and large $k$,  similar to  \eqref{eq:jnukasympt}, is well-known, see \cite[(10.21.20)]{DLMF}, we failed to find any analogues of \eqref{eq:jnukHEL} in this case. There is an analogue of the upper bound in \eqref{eq:jnukpQW}:
if
\begin{equation*}
(0>)a'_1>a'_2>\dots
\end{equation*}
denotes the  sequence of negative zeros of the derivative of the Airy function $\operatorname{Ai}'(x)$, then 
\begin{equation*}
\begin{split}
j'_{\nu,k}&<\nu + \left|a'_k\right| \left(\frac{\nu }{2}+\frac{8\left|a'_k\right|^{3/2}}{27}\right)^{1/3}+\frac{9\left|a'_k\right|^2}{10\times 2^{2/3}} \left(27 \nu +16 \left|a'_k\right|^{3/2}\right)^{-1/3}=:\wtilde{\ell}'_{\nu,k}\\
&\qquad\qquad\qquad\qquad\qquad\qquad\qquad\qquad\qquad\qquad\qquad\qquad\text{for  all }\nu\ge 0\text{ and }k\in\mathbb{N},
\end{split}
\end{equation*}
which is obtained from \cite{EL97}, \cite[\S1.7]{E01} after some manipulations required by the change of notation.

We compare these bounds to our bounds from Corollaries \ref{cor:jbounds}  and \ref{cor:jprimebounds} numerically. The results of the computations are collated in Tables \ref{table:1}--\ref{table:3} in Appendix \ref{sec:numerical}. We note that all the calculations have been performed with $16$ digit precision, and any differences between bounds not exceeding $10^{-16}$ are ignored. The inverse functions were evaluated using standard \texttt{Mathematica} routines with fixed precision.

The upper bound for $j_{\nu,k}$ obtained by using McMahon's expansions slightly outperforms our bound for low values of $\nu$ ($\lesssim 5$) and $k$  ($\lesssim 100$), however even in the worst case $\nu=\frac{1}{2}$, $k=1$, the relative deficiency of our bound is $7\times 10^{-4}$. For a fixed $\nu$ and very  large $k$, the bound $\wtilde{\ell}_{\nu,k}$ is theoretically better than our bound $\wtilde{j}_{\nu,k}$,  however in practice the difference is negligible. 
Similarly, the Airy-type upper bound $\wtilde{q}_{\nu,k}$ should be theoretically better than our bound once one looks deep in the transition region, but in practice within the range of parameters we checked,  it outperforms our bound only for very large $\nu$ and very small $k$ ($k \lesssim 10$ for $\nu=500000$), with the maximal relative deficiency of our bound being $3\times 10^{-4}$: in all other cases our bound is better.

Elbert--Laforgia's lower bound $\utilde{\ell}_{\nu,k}$ generally outperforms our bound $\utilde{j}_{\nu, k}$ (as it is based on an additional term of the asymptotic expansion) but not by much, and is only valid for very low $\nu$ ($\lesssim 1$).  The same observation as in the case of the upper bounds applies to Qu--Wang's lower bound: it becomes slightly more efficient than our bound only for very large $\nu$ and very small $k$. 

Concerning the bounds for $j'_{\nu,k}$, there is no alternative lower bound to compare with; the upper bound $\wtilde{\ell}'_{\nu,k}$ only  performs slightly better than our bound for very small $k$ and very large $\nu$.

Since accurately computing the large zeros of Bessel functions of large order and their derivatives is a non-trivial task, we do not compare our bounds to the actual values of these zeros; instead we plot, in Figures \ref{fig:errj} and \ref{fig:errjprime}, the \emph{maximal} possible errors by comparing our upper and lower bounds.

\phantomsection\section*{Acknowledgements}\addcontentsline{toc}{section}{Acknowledgements}
Research of NF was supported by the RSF grant  22-11-00092. Research of ML was partially supported by the EPSRC grants EP/W006898/1 and 
EP/V051881/1, and by the University of Reading RETF Open Fund. Research of IP was partially supported by NSERC. 

\phantomsection

\appendix

\section{Explicit expressions for some potentials}\label{sec:potentials}

We have 
\begin{equation}\label{eq:QthetaDown1}
\mathcal{V}_{\utilde{\theta}_\nu}(x)=\frac{q_\nu^{(1)}(x)}{q_\nu^{(2)}(x)},
\end{equation}
where 
\begin{equation}\label{eq:QthetaDown2}
\begin{split}
q_\nu^{(1)}(x):=&\hphantom{+\,}4096 x^{24}-2048 \left(24 \nu ^2-1\right) x^{22}+128 \left(2112 \nu ^4-128 \nu ^2+15\right) x^{20}\\
&-32 \left(28160 \nu ^6-1760 \nu ^4-584 \nu ^2-1\right) x^{18}\\
&+\left(2027520 \nu ^8-107520 \nu ^6-117376 \nu ^4+704 \nu ^2+1\right) x^{16}\\
&-16 \nu ^2 \left(202752 \nu ^8-7680 \nu ^6-13088 \nu ^4+223 \nu ^2-1\right) x^{14}\\
&+16\nu ^4 \left(236544 \nu ^8-5376 \nu ^6-5000 \nu ^4+585 \nu ^2+6\right)  x^{12}\\
&-16 \nu ^6 \left(202752 \nu ^8-2688 \nu ^6+11440 \nu ^4+935 \nu ^2-16\right) x^{10}\\
&+16\nu ^8 \left(126720 \nu ^8-1920 \nu ^6+15272 \nu ^4+823 \nu ^2+16\right)  x^8\\
&-128 \nu ^{12} \left(7040 \nu ^6-240 \nu ^4+808 \nu ^2+43\right) x^6\\
&+256 \left(1056 \nu ^6-80 \nu ^4+17 \nu ^2+3\right) \nu ^{14} x^4\\
&-1024 \nu ^{18} \left(48 \nu ^4-7 \nu ^2-5\right) x^2+1024 \left(4 \nu ^2-1\right) \nu ^{22},
\end{split}
\end{equation}
and
\begin{equation}\label{eq:QthetaDown3}
q_\nu^{(2)}(x):=64 x^2 \left(x^2-\nu ^2\right)^5  \left(8 x^6+\left(1-24 \nu ^2\right) x^4+4 \nu ^2\left(6 \nu ^2+1\right) x^2-8 \nu ^6\right)^2.
\end{equation}

Also,
\begin{equation}\label{eq:QphiUp}
\mathcal{V}_{\wtilde{\phi}_\nu}(x)=\frac{q^{(3)}_\nu(x)}{q^{(4)}_\nu(x)},
\end{equation}
where
\begin{equation}\label{eq:QphiUp3}
\begin{split}
q^{(3)}_\nu(x):=&\hphantom{+\,}4096 x^{24}-6144 \left(8 \nu ^2+1\right) x^{22}+128 \left(2112 \nu ^4+320 \nu ^2-9\right) x^{20}\\
&-32 \left(28160 \nu ^6+2848 \nu ^4+760 \nu ^2+27\right) x^{18}\\
&+3 \left(675840 \nu ^8-3072 \nu ^6+46208 \nu ^4+448 \nu ^2+27\right) x^{16}\\
&-16 \nu ^2 \left(202752 \nu ^8-29184 \nu ^6+16352 \nu ^4+575 \nu ^2-27\right) x^{14}\\
&+16 \nu ^4 \left(236544 \nu ^8-69888 \nu ^6+10360 \nu ^4+1657 \nu ^2+54\right)  x^{12}\\
&-48 \nu ^6 \left(67584 \nu ^8-29568 \nu ^6-1904 \nu ^4+533 \nu ^2-16\right) x^{10}\\
&+16  \nu ^8 \left(126720 \nu ^8-69504 \nu ^6-11480 \nu ^4+479 \nu ^2+16\right) x^8\\
&-128 \nu ^{12} \left(7040 \nu ^6-4272 \nu ^4-632 \nu ^2+5\right) x^6\\
&+768 \nu ^{14} \left(352 \nu ^6-208 \nu ^4-\nu ^2+1\right)  x^4\\
&-1024 \nu ^{18} \left(48 \nu ^4-23 \nu ^2+5\right) x^2 + 1024 \left(4 \nu ^2-1\right) \nu ^{22}
\end{split}
\end{equation}
and
\begin{equation}\label{eq:QphiUp4}
q^{(4)}_\nu(x):=64 x^2 \left(x^2-\nu ^2\right)^5  \left(8 x^6-3 \left(8 \nu ^2+1\right) x^4+4 \left(6 \nu ^2-1\right) \nu ^2 x^2-8 \nu ^6\right)^2.
\end{equation}

Finally,
\begin{equation}\label{eq:QpsiDown}
\mathcal{V}_{\utilde{\psi}_{\nu,\eta}}(x)=\frac{q^{(5)}_{\mu,\eta}(x)}{q^{(6)}_{\mu,\eta}(x)},\qquad \mu=\mu_{\nu,\eta},
\end{equation}
where
\begin{equation}\label{eq:QpsiDown5}
\begin{split}
q^{(5)}_{\mu,\eta}(x):=&\hphantom{+\,}16 x^{16}\\
&-32 x^{14} \left(\eta ^2+2 \eta +4 \mu ^2\right)\\
&+8 x^{12} \left(3 \eta ^4+12 \eta ^3+\eta ^2 \left(28 \mu ^2+9\right)+6 \eta  \left(8 \mu ^2-1\right)+\left(56 \mu ^2-3\right) \mu ^2\right)\\
&-4 x^{10} \left(2 \eta ^6+12 \eta ^5+12 \eta ^4 \left(3 \mu ^2+2\right)+8 \eta ^3 \left(15 \mu ^2+2\right)\right.\\
&\left.\qquad+14 \eta ^2 \left(12 \mu ^2+5\right) \mu ^2+12 \eta  \left(20 \mu ^2-1\right) \mu ^2+\left(224 \mu ^2-31\right) \mu ^4\right)\\
&+x^8 \left(\eta ^8+8 \eta ^7+8 \eta ^6 \left(5 \mu ^2+3\right)+32 \eta ^5 \left(6 \mu ^2+1\right)+2 \eta ^4 \left(180 \mu ^4+145 \mu ^2+8\right)\right.\\
&\left.\qquad+48 \eta ^3 \left(20 \mu ^2+3\right) \mu ^2+4 \eta ^2 \left(280 \mu ^4+99 \mu ^2+6\right) \mu ^2\right.\\
&\left.\qquad+8 \eta  \left(160 \mu ^2+13\right) \mu ^4+20 \left(56 \mu ^2-13\right) \mu ^6\right)\\
&-\mu ^2 x^6 \left(4 \eta ^8+24 \eta ^7+16 \eta ^6 \left(5 \mu ^2+3\right)+32 \eta ^5 \left(9 \mu ^2+1\right)+\eta ^4 \left(480 \mu ^2+293\right) \mu ^2\right.\\
&\left.\qquad+\eta ^3 \left(960 \mu ^4+76 \mu ^2\right)+4 \eta ^2 \left(280 \mu ^4+56 \mu ^2+9\right) \mu ^2\right.\\
&\left.\qquad+120 \eta  \left(8 \mu ^6+\mu ^4\right)+56 \left(16 \mu ^2-5\right) \mu ^6\right)\\
&+\mu ^4 x^4 \left(6 \eta ^8+24 \eta ^7+8 \eta ^6 \left(10 \mu ^2+3\right)+192 \eta ^5 \mu ^2+9 \eta ^4 \left(40 \mu ^2+11\right) \mu ^2\right.\\
&\left.\qquad+24 \eta ^3 \left(20 \mu ^2-1\right) \mu ^2+4 \eta ^2 \left(168 \mu ^4+4 \mu ^2+3\right) \mu ^2\right.\\
&\left.\qquad+24 \eta  \left(16 \mu ^2-1\right) \mu ^4+32 \left(14 \mu ^2-5\right) \mu ^6\right)\\
&-\mu ^6 x^2 \left(4 \eta ^8+8 \eta ^7+40 \eta ^6 \mu ^2+48 \eta ^5 \mu ^2+\eta ^4 \left(144 \mu ^2-1\right) \mu ^2+4 \eta ^3 \left(24 \mu ^2-5\right) \mu ^2\right.\\
&\left.\qquad+8 \eta ^2 \left(28 \mu ^2-3\right) \mu ^4+8 \eta  \left(8 \mu ^2-5\right) \mu ^4+4 \left(32 \mu ^2-11\right) \mu ^6\right)\\
&+\mu ^8 \left(\eta ^4+4 \eta ^2 \mu ^2+4 \mu ^4-\mu ^2\right) \left(\eta ^2+2 \mu ^2\right)^2
\end{split}
\end{equation}
and
\begin{equation}\label{eq:QpsiDown6}
q^{(6)}_{\mu,\eta}(x):=4 (x^2-\mu^2 )^3 x^2 \left(2 x^4-x^2 \left((\eta +2) \eta +4 \mu ^2\right)+\mu ^2 \left(\eta ^2+2 \mu ^2\right)\right)^2.
\end{equation}

\section{Numerical data}\label{sec:numerical}

\subsection{Bounds for zeros of Bessel functions}
{\small\begin{longtable}{@{}>{$}r<{$}>{$}r<{$}@{\hspace{1cm}}*{4}{>{$}c<{$}}@{}}
\caption{The comparison of our, Elbert--Laforgia's, and Qu--Wong's upper bounds for $j_{\nu, k}$. If $\wtilde{m}_{\nu,k}:=\min\left\{\wtilde{j}_{\nu,k}, \wtilde{\ell}_{\nu,k}, \wtilde{q}_{\nu,k}\right\}$ denotes the best of the three bounds, the coloured entries in each row are those matching it. If $\wtilde{j}_{\nu,k}>\wtilde{m}_{\nu,k}$, the entry in the last column shows the relative deficiency of our bound.\label{table:1}}\\[1cm]
\nopagebreak
\toprule
\nu & k & \wtilde{j}_{\nu,k} & \wtilde{\ell}_{\nu,k} & \wtilde{q}_{\nu,k} & \frac{\wtilde{j}_{\nu,k}}{\wtilde{m}_{\nu,k}}-1\\\nopagebreak\midrule
\endfirsthead
\toprule
\nu & k & \wtilde{j}_{\nu,k} & \wtilde{\ell}_{\nu,k} & \wtilde{q}_{\nu,k} &  \frac{\wtilde{j}_{\nu,k}}{\wtilde{m}_{\nu,k}}-1 \\\nopagebreak\midrule
\endhead
\multirow{10}{*}{$0$} & 1 & {\color{red}2.40810257797209} & 2.40924613788964 &   & \\\nopagebreak
 & 2 & {\color{red}5.52043030619282} & 5.52052356422384 &   & \\\nopagebreak
 & 5 & {\color{red}14.9309369836999} & 14.9309416805037 &   & \\\nopagebreak
 & 10 & {\color{red}30.6346087249255} & 30.6346092684772 &   & \\\nopagebreak
 & 50 & {\color{red}156.295034285583} & 156.295034289676 &   & \\\nopagebreak
 & 100 & {\color{red}313.374266079643} & 313.374266080151 &   & \\\nopagebreak
 & 1000 & {\color{red}3140.80729522508} & {\color{red}3140.80729522508} &   & \\\nopagebreak
 & 10000 & {\color{red}31415.1411417135} & {\color{red}31415.1411417135} &   & \\\nopagebreak
 & 100000 & {\color{red}314158.479961214} & {\color{red}314158.479961214} &   & \\\nopagebreak
 & 500000 & {\color{red}1570795.54139681} & {\color{red}1570795.54139681} &   & \\
\midrule
\multirow{10}{*}{$\frac{1}{2}$} & 1 & 3.14376279614743 & {\color{red}3.14159265358979} & 3.27460318694968 & 7\times 10^{-4}\\\nopagebreak
 & 2 & 6.28345002833208 & {\color{red}6.28318530717959} & 7.05438435767338 & 4\times 10^{-5}\\\nopagebreak
 & 5 & 15.7079800888803 & {\color{red}15.7079632679490} & 20.5314460526698 & 1\times 10^{-6}\\\nopagebreak
 & 10 & 31.4159286363346 & {\color{red}31.4159265358979} & 47.7692005216613 & 7\times 10^{-8}\\\nopagebreak
 & 50 & 157.079632696288 & {\color{red}157.079632679490} & 368.663088373722 & 1\times 10^{-10}\\\nopagebreak
 & 100 & 314.159265361079 & {\color{red}314.159265358979} & 908.847335969838 & 7\times 10^{-12}\\\nopagebreak
 & 1000 & {\color{red}3141.59265358980} & {\color{red}3141.59265358979} & 18983.1730584802 & \\\nopagebreak
 & 10000 & {\color{red}31415.9265358979} & {\color{red}31415.9265358979} & 406099.051081796 & \\\nopagebreak
 & 100000 & {\color{red}314159.265358979} & {\color{red}314159.265358979} & 8735498.79115534 & \\\nopagebreak
 & 500000 & {\color{red}1570796.32679490} & {\color{red}1570796.32679490} & 74666197.7784502 & \\
\midrule
\multirow{10}{*}{$1$} & 1 & 3.83369791861688 & {\color{red}3.83188486954530} & 3.88890738513848 & 5\times 10^{-4}\\\nopagebreak
 & 2 & 7.01581736210639 & {\color{red}7.01559818406170} & 7.40285121412499 & 3\times 10^{-5}\\\nopagebreak
 & 5 & 16.4706452289041 & {\color{red}16.4706302346712} & 19.2321654812515 & 9\times 10^{-7}\\\nopagebreak
 & 10 & 32.1896818843079 & {\color{red}32.1896799175785} & 42.2854077969397 & 6\times 10^{-8}\\\nopagebreak
 & 50 & 157.862655418487 & {\color{red}157.862655401933} & 304.378027137110 & 1\times 10^{-10}\\\nopagebreak
 & 100 & 314.943472839851 & {\color{red}314.943472837767} & 739.711587305586 & 7\times 10^{-12}\\\nopagebreak
 & 1000 & {\color{red}3142.37793241682} & {\color{red}3142.37793241682} & 15150.0966962282 & \\\nopagebreak
 & 10000 & {\color{red}31416.7119221250} & {\color{red}31416.7119221250} & 322704.803661024 & \\\nopagebreak
 & 100000 & {\color{red}314160.050755949} & {\color{red}314160.050755949} & 6935149.13300128 & \\\nopagebreak
 & 500000 & {\color{red}1570797.11219282} & {\color{red}1570797.11219282} & 59267801.5559670 & \\
\midrule
\multirow{10}{*}{$5$} & 1 & {\color{red}8.77372338204280} & 8.83610697727954 & 8.77748993504362 & \\\nopagebreak
 & 2 & {\color{red}12.3388145653858} & 12.3524627245014 & 12.3951530492719 & \\\nopagebreak
 & 5 & {\color{red}22.2178113198726} & 22.2186275958908 & 22.7567429583036 & \\\nopagebreak
 & 10 & {\color{red}38.1598700502104} & 38.1599260373733 & 40.6005910616904 & \\\nopagebreak
 & 50 & {\color{red}164.072787945468} & 164.072787970360 & 216.371205491807 & \\\nopagebreak
 & 100 & 321.189319569572 & {\color{red}321.189319568987} & 490.990987826608 & 2\times 10^{-12}\\\nopagebreak
 & 1000 & {\color{red}3148.65730681305} & {\color{red}3148.65730681305} & 9115.22066902932 & \\\nopagebreak
 & 10000 & {\color{red}31422.9947255486} & {\color{red}31422.9947255486} & 189888.427570357 & \\\nopagebreak
 & 100000 & {\color{red}314166.333903060} & {\color{red}314166.333903060} & 4061112.22022502 & \\\nopagebreak
 & 500000 & {\color{red}1570803.39537049} & {\color{red}1570803.39537049} & 34675838.1018364 & \\
\midrule
\multirow{10}{*}{$10$} & 1 & 14.4780593618721 & 14.8150624867862 & {\color{red}14.4776533239189} & 3\times 10^{-5}\\\nopagebreak
 & 2 & {\color{red}18.4336923418614} & 18.5520743143216 & 18.4562220388682 & \\\nopagebreak
 & 5 & {\color{red}28.8873862860562} & 28.9021354143944 & 29.1202550857628 & \\\nopagebreak
 & 10 & {\color{red}45.2315754422830} & 45.2332565697875 & 46.3737266901699 & \\\nopagebreak
 & 50 & {\color{red}171.711662928231} & 171.711665150058 & 201.823557261092 & \\\nopagebreak
 & 100 & {\color{red}328.930191596730} & 328.930191681758 & 433.985076544903 & \\\nopagebreak
 & 1000 & {\color{red}3156.49941795039} & {\color{red}3156.49941795039} & 7418.61204821241 & \\\nopagebreak
 & 10000 & {\color{red}31430.8475141856} & {\color{red}31430.8475141856} & 151546.091220676 & \\\nopagebreak
 & 100000 & {\color{red}314174.187765333} & {\color{red}314174.187765334} & 3227144.69490720 & \\\nopagebreak
 & 500000 & {\color{red}1570811.24932825} & {\color{red}1570811.24932825} & 27533441.1155637 & \\
\midrule
\multirow{10}{*}{$50$} & 1 & 57.1207457655499 & 62.0041045949062 & {\color{red}57.1171076805129} & 6\times 10^{-5}\\\nopagebreak
 & 2 & {\color{red}62.8080187024269} & 66.0949975045479 & 62.8105158914088 & \\\nopagebreak
 & 5 & {\color{red}76.4370860305146} & 77.8577326931975 & 76.4662469423382 & \\\nopagebreak
 & 10 & {\color{red}95.8011096962304} & 96.3212070986120 & 95.9542784414793 & \\\nopagebreak
 & 50 & {\color{red}229.362879678685} & 229.370990348625 & 235.332191927489 & \\\nopagebreak
 & 100 & {\color{red}388.693660067007} & 388.694258337930 & 414.265541479911 & \\\nopagebreak
 & 1000 & {\color{red}3218.95877848403} & 3218.95877849964 & 4923.29223763353 & \\\nopagebreak
 & 10000 & {\color{red}31493.6412674832} & {\color{red}31493.6412674832} & 91178.9720478066 & \\\nopagebreak
 & 100000 & {\color{red}314237.015799664} & {\color{red}314237.015799664} & 1898940.97658646 & \\\nopagebreak
 & 500000 & {\color{red}1570874.08041742} & {\color{red}1570874.08041742} & 16135766.4340724 & \\
\midrule
\multirow{10}{*}{$100$} & 1 & 108.840879423458 & 120.880380157578 & {\color{red}108.836246832758} & 4\times 10^{-5}\\\nopagebreak
 & 2 & {\color{red}115.739736248993} & 125.037089659834 & 115.740557472041 & \\\nopagebreak
 & 5 & {\color{red}131.823950776853} & 137.202588450753 & 131.836011607965 & \\\nopagebreak
 & 10 & {\color{red}153.900272851252} & 156.664768595604 & 153.962599154298 & \\\nopagebreak
 & 50 & {\color{red}296.335776171518} & 296.471251735914 & 298.929907132201 & \\\nopagebreak
 & 100 & {\color{red}459.529546576725} & 459.545645439361 & 471.533009349504 & \\\nopagebreak
 & 1000 & {\color{red}3296.36998972096} & 3296.36999060841 & 4351.04720098391 & \\\nopagebreak
 & 10000 & {\color{red}31572.0624063796} & {\color{red}31572.0624063796} & 74207.6309465587 & \\\nopagebreak
 & 100000 & {\color{red}314315.543686312} & {\color{red}314315.543686312} & 1515506.03955682 & \\\nopagebreak
 & 500000 & {\color{red}1570952.61784671} & {\color{red}1570952.61784671} & 12831165.2167116 & \\
\midrule
\multirow{10}{*}{$1000$} & 1 & 1018.67064844448 & 1180.40361331412 & {\color{red}1018.66088584526} & 1\times 10^{-5}\\\nopagebreak
 & 2 & 1032.76258530442 & 1184.62569023296 & {\color{red}1032.76190059904} & 7\times 10^{-7}\\\nopagebreak
 & 5 & {\color{red}1064.24453144028} & 1197.25570362592 & 1064.24532031332 & \\\nopagebreak
 & 10 & {\color{red}1104.92859950901} & 1218.18794810042 & 1104.93238865887 & \\\nopagebreak
 & 50 & {\color{red}1328.95755863688} & 1380.97019620599 & 1329.09301911099 & \\\nopagebreak
 & 100 & {\color{red}1548.25088464171} & 1575.19750272123 & 1548.90887404113 & \\\nopagebreak
 & 1000 & {\color{red}4602.53426352435} & 4602.69406661919 & 4723.15283247506 & \\\nopagebreak
 & 10000 & {\color{red}32970.7713584907} & 32970.7713673585 & 43521.6749153743 & \\\nopagebreak
 & 100000 & {\color{red}315727.692643526} & {\color{red}315727.692643526} & 742097.821357475 & \\\nopagebreak
 & 500000 & {\color{red}1572366.01973148} & {\color{red}1572366.01973148} & 6066913.80070566 & \\
\midrule
\multirow{10}{*}{$10000$} & 1 & 10040.0498904506 & 11775.5007480898 & {\color{red}10040.0290289155} & 2\times 10^{-6}\\\nopagebreak
 & 2 & 10070.0511840226 & 11779.7297125815 & {\color{red}10070.0495448888} & 2\times 10^{-7}\\\nopagebreak
 & 5 & {\color{red}10136.3963629401} & 11792.4129183083 & 10136.3963712851 & \\\nopagebreak
 & 10 & {\color{red}10220.8123677153} & 11813.5393327241 & 10220.8126794286 & \\\nopagebreak
 & 50 & {\color{red}10662.8217778127} & 11982.0061387084 & 10662.8309577604 & \\\nopagebreak
 & 100 & {\color{red}11065.7962862662} & 12191.2634222777 & 11065.8362405931 & \\\nopagebreak
 & 1000 & {\color{red}15491.7628369823} & 15760.5279241169 & 15498.3768788958 & \\\nopagebreak
 & 10000 & {\color{red}46032.5838776516} & 46034.1806805361 & 47239.3537167537 & \\\nopagebreak
 & 100000 & {\color{red}329714.785384131} & 329714.785472800 & 435227.952705537 & \\\nopagebreak
 & 500000 & {\color{red}1586471.98808840} & 1586471.98808844 & 3063097.67857318 & \\
\midrule
\multirow{10}{*}{$100000$} & 1 & 100086.203733327 & 117726.458311625 & {\color{red}100086.158872022} & 4\times 10^{-7}\\\nopagebreak
 & 2 & 100150.672924750 & 117730.687968554 & {\color{red}100150.669387466} & 4\times 10^{-8}\\\nopagebreak
 & 5 & 100292.921471816 & 117743.376569896 & {\color{red}100292.921340968} & 1\times 10^{-9}\\\nopagebreak
 & 10 & {\color{red}100473.286256924} & 117764.523007627 & 100473.286274617 & \\\nopagebreak
 & 50 & {\color{red}101406.591053031} & 117933.639180549 & 101406.591867126 & \\\nopagebreak
 & 100 & {\color{red}102242.079713765} & 118144.896439004 & 102242.083055865 & \\\nopagebreak
 & 1000 & {\color{red}110674.460043079} & 121922.018138173 & 110674.861526695 & \\\nopagebreak
 & 10000 & {\color{red}154926.882893420} & 157613.832441531 & 154993.057404632 & \\\nopagebreak
 & 100000 & {\color{red}460333.080262841} & 460349.047059410 & 472401.362824375 & \\\nopagebreak
 & 500000 & {\color{red}1724975.77003468} & 1724975.79260235 & 2028908.15934764 & \\
\midrule
\multirow{10}{*}{$500000$} & 1 & 500147.381239784 & 588619.602430817 & {\color{red}500147.304554046} & 2\times 10^{-7}\\\nopagebreak
 & 2 & 500257.570515678 & 588623.832149329 & {\color{red}500257.564469154} & 1\times 10^{-8}\\\nopagebreak
 & 5 & 500500.599560184 & 588636.521230963 & {\color{red}500500.599326034} & 5\times 10^{-10}\\\nopagebreak
 & 10 & 500808.554185199 & 588657.669454029 & {\color{red}500808.554169543} & 3\times 10^{-11}\\\nopagebreak
 & 50 & {\color{red}502398.615413161} & 588826.844156947 & 502398.615572073 & \\\nopagebreak
 & 100 & {\color{red}503817.163442824} & 589038.284846206 & 503817.164086853 & \\\nopagebreak
 & 1000 & {\color{red}517891.863383625} & 592838.985093676 & 517891.932705137 & \\\nopagebreak
 & 10000 & {\color{red}586230.238693899} & 630329.382318621 & 586239.210056145 & \\\nopagebreak
 & 100000 & {\color{red}967172.336060185} & 972162.043475722 & 968798.233994082 & \\\nopagebreak
 & 500000 & {\color{red}2301668.61976330} & 2301748.45319828 & 2362010.29220323 & \\
\bottomrule
\end{longtable}}

{\small\begin{longtable}{@{}>{$}r<{$}>{$}r<{$}@{\hspace{1cm}}*{4}{>{$}c<{$}}@{}}
\caption{The comparison of our, Elbert--Laforgia's, and Qu--Wong's lower bounds for $j_{\nu, k}$. If $\utilde{m}_{\nu,k}:=\max\left\{\utilde{j}_{\nu,k}, \utilde{\ell}_{\nu,k}, \utilde{q}_{\nu,k}\right\}$ denotes the best of the three bounds, the coloured entries in each row are those matching it. If $\utilde{j}_{\nu,k}<\utilde{m}_{\nu,k}$, the entry in the last column shows the relative deficiency of our bound.\label{table:2}}\\[1cm]
\nopagebreak
\toprule
\nu & k & \utilde{j}_{\nu,k} & \utilde{\ell}_{\nu,k} & \utilde{q}_{\nu,k} &  \frac{\utilde{j}_{\nu,k}}{\utilde{m}_{\nu,k}}-1 \\\nopagebreak\midrule
\endfirsthead
\toprule
\nu & k & \utilde{j}_{\nu,k} & \utilde{\ell}_{\nu,k} & \utilde{q}_{\nu,k} &  \frac{\utilde{j}_{\nu,k}}{\utilde{m}_{\nu,k}}-1 \\\nopagebreak\midrule
\endhead
\multirow{10}{*}{$0$} & 1 & 2.35619449019234 & {\color{red}2.40307454796724} &   & -2\times 10^{-2}\\\nopagebreak
 & 2 & 5.49778714378214 & {\color{red}5.52003775393840} &   & -4\times 10^{-3}\\\nopagebreak
 & 5 & 14.9225651045515 & {\color{red}14.9309173864480} &   & -6\times 10^{-4}\\\nopagebreak
 & 10 & 30.6305283725005 & {\color{red}30.6346064593785} &   & -1\times 10^{-4}\\\nopagebreak
 & 50 & 156.294234516092 & {\color{red}156.295034268531} &   & -5\times 10^{-6}\\\nopagebreak
 & 100 & 313.373867195582 & {\color{red}313.374266077528} &   & -1\times 10^{-6}\\\nopagebreak
 & 1000 & 3140.80725542640 & {\color{red}3140.80729522508} &   & -1\times 10^{-8}\\\nopagebreak
 & 10000 & 31415.1411377345 & {\color{red}31415.1411417135} &   & -1\times 10^{-10}\\\nopagebreak
 & 100000 & 314158.479960816 & {\color{red}314158.479961214} &   & -1\times 10^{-12}\\\nopagebreak
 & 500000 & 1570795.54139673 & {\color{red}1570795.54139681} &   & -5\times 10^{-14}\\
\midrule
\multirow{10}{*}{$\frac{1}{2}$} & 1 & 3.10119764278657 & {\color{red}3.14159265358979} & 1.97291537167673 & -1\times 10^{-2}\\\nopagebreak
 & 2 & 6.26321689238053 & {\color{red}6.28318530717959} & 3.07524677778333 & -3\times 10^{-3}\\\nopagebreak
 & 5 & 15.7000008117865 & {\color{red}15.7079632679490} & 5.50449056479522 & -5\times 10^{-4}\\\nopagebreak
 & 10 & 31.4119470742254 & {\color{red}31.4159265358979} & 8.58162293766878 & -1\times 10^{-4}\\\nopagebreak
 & 50 & 157.078836900071 & {\color{red}157.079632679490} & 24.4517345853548 & -5\times 10^{-6}\\\nopagebreak
 & 100 & 314.158867471034 & {\color{red}314.159265358979} & 38.5846145963924 & -1\times 10^{-6}\\\nopagebreak
 & 1000 & 3141.59261380106 & {\color{red}3141.59265358979} & 177.538763621880 & -1\times 10^{-8}\\\nopagebreak
 & 10000 & 31415.9265319191 & {\color{red}31415.9265358979} & 822.364431896465 & -1\times 10^{-10}\\\nopagebreak
 & 100000 & 314159.265358581 & {\color{red}314159.265358979} & 3815.31399222450 & -1\times 10^{-12}\\\nopagebreak
 & 500000 & 1570796.32679482 & {\color{red}1570796.32679490} & 11155.0986540527 & -5\times 10^{-14}\\
\midrule
\multirow{10}{*}{$1$} & 1 & 3.79443997608576 & {\color{red}3.83149785113210} & 2.85575708148924 & -1\times 10^{-2}\\\nopagebreak
 & 2 & 6.99700190767405 & {\color{red}7.01553182287974} & 4.24460762400316 & -3\times 10^{-3}\\\nopagebreak
 & 5 & 16.4629809123032 & {\color{red}16.4706250109047} & 7.30526300658577 & -5\times 10^{-4}\\\nopagebreak
 & 10 & 32.1857886427175 & {\color{red}32.1896792156545} & 11.1822068564821 & -1\times 10^{-4}\\\nopagebreak
 & 50 & 157.861863506199 & {\color{red}157.862655395975} & 31.1772945855835 & -5\times 10^{-6}\\\nopagebreak
 & 100 & 314.943075932560 & {\color{red}314.943472837017} & 48.9836076071283 & -1\times 10^{-6}\\\nopagebreak
 & 1000 & 3142.37789263802 & {\color{red}3142.37793241682} & 224.054864934569 & -1\times 10^{-8}\\\nopagebreak
 & 10000 & 31416.7119181462 & {\color{red}31416.7119221250} & 1036.48429790625 & -1\times 10^{-10}\\\nopagebreak
 & 100000 & 314160.050755551 & {\color{red}314160.050755949} & 4807.36445023715 & -1\times 10^{-12}\\\nopagebreak
 & 500000 & 1570797.11219274 & {\color{red}1570797.11219282} & 14054.9136473700 & -5\times 10^{-14}\\
\midrule
\multirow{10}{*}{$5$} & 1 & {\color{red}8.73567022436826} &   & 8.17329997222154 & \\\nopagebreak
 & 2 & {\color{red}12.3227225053509} &   & 10.5482009934492 & \\\nopagebreak
 & 5 & {\color{red}22.2113581058685} &   & 15.7818480787321 & \\\nopagebreak
 & 10 & {\color{red}38.1564401016457} &   & 22.4113288086710 & \\\nopagebreak
 & 50 & {\color{red}164.072024198901} &   & 56.6024478771247 & \\\nopagebreak
 & 100 & {\color{red}321.188930139405} &   & 87.0508148429624 & \\\nopagebreak
 & 1000 & {\color{red}3148.65726711332} &   & 386.418453827333 & \\\nopagebreak
 & 10000 & {\color{red}31422.9947215707} &   & 1775.65324258109 & \\\nopagebreak
 & 100000 & {\color{red}314166.333902662} &   & 8223.76760086749 & \\\nopagebreak
 & 500000 & {\color{red}1570803.39537041} &   & 24036.8542936740 & \\
\midrule
\multirow{10}{*}{$10$} & 1 & {\color{red}14.4363906385824} &   & 13.9981074326327 & \\\nopagebreak
 & 2 & {\color{red}18.4174021405972} &   & 16.9902952206943 & \\\nopagebreak
 & 5 & {\color{red}28.8813543155782} &   & 23.5842773511631 & \\\nopagebreak
 & 10 & {\color{red}45.2284209970154} &   & 31.9368996726856 & \\\nopagebreak
 & 50 & {\color{red}171.710928343398} &   & 75.0150103064925 & \\\nopagebreak
 & 100 & {\color{red}328.929810638849} &   & 113.377548781675 & \\\nopagebreak
 & 1000 & {\color{red}3156.49937834850} &   & 490.557138795412 & \\\nopagebreak
 & 10000 & {\color{red}31430.8475102086} &   & 2240.88329239253 & \\\nopagebreak
 & 100000 & {\color{red}314174.187764936} &   & 10364.9983045269 & \\\nopagebreak
 & 500000 & {\color{red}1570811.24932817} &   & 30288.2390926064 & \\
\midrule
\multirow{10}{*}{$50$} & 1 & {\color{red}57.0601473524213} &   & 56.8366675420313 & \\\nopagebreak
 & 2 & {\color{red}62.7869087299433} &   & 61.9532366875564 & \\\nopagebreak
 & 5 & {\color{red}76.4306643187313} &   & 73.2287875234740 & \\\nopagebreak
 & 10 & {\color{red}95.7981975064759} &   & 87.5115707849523 & \\\nopagebreak
 & 50 & {\color{red}229.362259920245} &   & 161.174103797040 & \\\nopagebreak
 & 100 & {\color{red}388.693323896030} &   & 226.773121843061 & \\\nopagebreak
 & 1000 & {\color{red}3218.95873962660} &   & 871.741148343930 & \\\nopagebreak
 & 10000 & {\color{red}31493.6412635141} &   & 3864.75676983414 & \\\nopagebreak
 & 100000 & {\color{red}314237.015799267} &   & 17756.7980286190 & \\\nopagebreak
 & 500000 & {\color{red}1570874.08041734} &   & 51825.0605560830 & \\
\midrule
\multirow{10}{*}{$100$} & 1 & {\color{red}108.766857714791} &   & 108.613661347338 & \\\nopagebreak
 & 2 & {\color{red}115.714556094408} &   & 115.060134517028 & \\\nopagebreak
 & 5 & {\color{red}131.816671193838} &   & 129.266438364360 & \\\nopagebreak
 & 10 & {\color{red}153.897155174351} &   & 147.261617646583 & \\\nopagebreak
 & 50 & {\color{red}296.335198187589} &   & 240.070593577088 & \\\nopagebreak
 & 100 & {\color{red}459.529237380598} &   & 322.720177265704 & \\\nopagebreak
 & 1000 & {\color{red}3296.36995170726} &   & 1135.32897036330 & \\\nopagebreak
 & 10000 & {\color{red}31572.0624024203} &   & 4906.29235454301 & \\\nopagebreak
 & 100000 & {\color{red}314315.543685914} &   & 22409.1675624942 & \\\nopagebreak
 & 500000 & {\color{red}1570952.61784663} &   & 65332.4886541907 & \\
\midrule
\multirow{10}{*}{$1000$} & 1 & 1018.51756702735 &   & {\color{red}1018.55757081489} & -4\times 10^{-5}\\\nopagebreak
 & 2 & {\color{red}1032.71211855974} &   & 1032.44607624003 & \\\nopagebreak
 & 5 & {\color{red}1064.23092084051} &   & 1063.05263006586 & \\\nopagebreak
 & 10 & {\color{red}1104.92325532692} &   & 1101.82206856482 & \\\nopagebreak
 & 50 & {\color{red}1328.95687010529} &   & 1301.77294585584 & \\\nopagebreak
 & 100 & {\color{red}1548.25058086023} &   & 1479.83607607128 & \\\nopagebreak
 & 1000 & {\color{red}4602.53423266599} &   & 3230.54864934569 & \\\nopagebreak
 & 10000 & {\color{red}32970.7713546901} &   & 11354.8429790625 & \\\nopagebreak
 & 100000 & {\color{red}315727.692643130} &   & 49063.6445023715 & \\\nopagebreak
 & 500000 & {\color{red}1572366.01973140} &   & 141539.136473700 & \\
\midrule
\multirow{10}{*}{$10000$} & 1 & 10039.7229535966 &   & {\color{red}10039.9810743263} & -3\times 10^{-5}\\\nopagebreak
 & 2 & {\color{red}10069.9441122528} &   & 10069.9029522069 & \\\nopagebreak
 & 5 & {\color{red}10136.3679056193} &   & 10135.8427735116 & \\\nopagebreak
 & 10 & {\color{red}10220.8014007353} &   & 10219.3689967269 & \\\nopagebreak
 & 50 & {\color{red}10662.8204962035} &   & 10650.1501030649 & \\\nopagebreak
 & 100 & {\color{red}11065.7957673197} &   & 11033.7754878167 & \\\nopagebreak
 & 1000 & {\color{red}15491.7628066820} &   & 14805.5713879541 & \\\nopagebreak
 & 10000 & {\color{red}46032.5838745664} &   & 32308.8329239253 & \\\nopagebreak
 & 100000 & {\color{red}329714.785383751} &   & 113549.983045269 & \\\nopagebreak
 & 500000 & {\color{red}1586471.98808833} &   & 312782.390926064 & \\
\midrule
\multirow{10}{*}{$100000$} & 1 & 100085.500689806 &   & {\color{red}100086.136613473} & -6\times 10^{-6}\\\nopagebreak
 & 2 & 100150.443003779 &   & {\color{red}100150.601345170} & -2\times 10^{-6}\\\nopagebreak
 & 5 & {\color{red}100292.860554032} &   & 100292.664383644 & \\\nopagebreak
 & 10 & {\color{red}100473.262872899} &   & 100472.616176466 & \\\nopagebreak
 & 50 & {\color{red}101406.588375922} &   & 101400.705935771 & \\\nopagebreak
 & 100 & {\color{red}102242.078649619} &   & 102227.201772657 & \\\nopagebreak
 & 1000 & {\color{red}110674.459991335} &   & 110353.289703633 & \\\nopagebreak
 & 10000 & {\color{red}154926.882890391} &   & 148062.923545430 & \\\nopagebreak
 & 100000 & {\color{red}460333.080262533} &   & 323091.675624942 & \\\nopagebreak
 & 500000 & {\color{red}1724975.77003461} &   & 752324.886541907 & \\
\midrule
\multirow{10}{*}{$500000$} & 1 & 500146.179459647 &   & {\color{red}500147.291537168} & -2\times 10^{-6}\\\nopagebreak
 & 2 & 500257.177589753 &   & {\color{red}500257.524677778} & -7\times 10^{-7}\\\nopagebreak
 & 5 & {\color{red}500500.495512600} &   & 500500.449056480 & \\\nopagebreak
 & 10 & {\color{red}500808.514273685} &   & 500808.162293767 & \\\nopagebreak
 & 50 & {\color{red}502398.610860722} &   & 502395.173458535 & \\\nopagebreak
 & 100 & {\color{red}503817.161639156} &   & 503808.461459639 & \\\nopagebreak
 & 1000 & {\color{red}517891.863298728} &   & 517703.876362188 & \\\nopagebreak
 & 10000 & {\color{red}586230.238689636} &   & 582186.443189646 & \\\nopagebreak
 & 100000 & {\color{red}967172.336059901} &   & 881481.399222450 & \\\nopagebreak
 & 500000 & {\color{red}2301668.61976324} &   & 1615459.86540527 & \\
\bottomrule
\end{longtable}}

\begin{figure}[ht]
\centering
\includegraphics{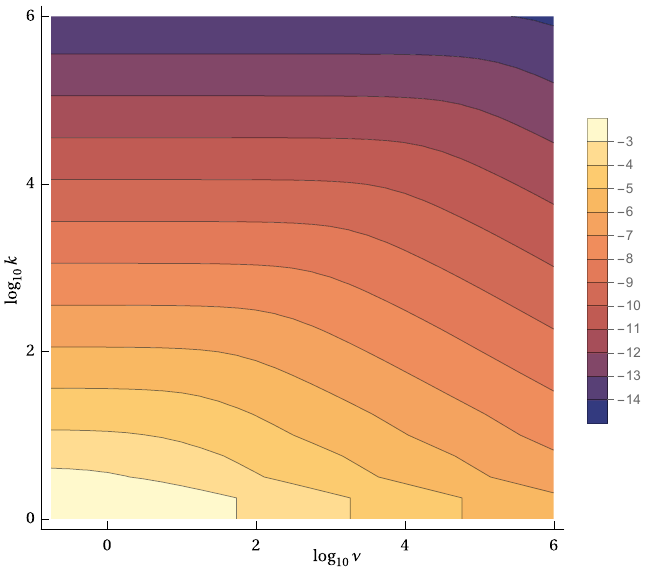}
\caption{The contour plot of $\log_{10}\left(\frac{\wtilde{j}_{\nu,k} }{\utilde{j}_{\nu,k}}-1\right)$, which estimates from above the maximum possible order of relative error of our approximations of $j_{\nu,k}$.\label{fig:errj}}
\end{figure}

\subsection{Bounds for zeros of derivatives of Bessel functions}

{\small\begin{longtable}{@{}>{$}r<{$}>{$}r<{$}@{\hspace{1cm}}>{$}c<{$}@{\hspace{1cm}}*{3}{>{$}c<{$}}@{}}
\caption{Our lower bound and the comparison of our and Elbert--Laforgia's upper bounds for $j'_{\nu, k}$. If $\wtilde{m}'_{\nu,k}:=\min\left\{\wtilde{j}'_{\nu,k}, \wtilde{\ell}'_{\nu,k}\right\}$ denotes the best of the two upper bounds, the coloured entries in each row are those matching it. If $\wtilde{j}'_{\nu,k}>\wtilde{m}'_{\nu,k}$, the entry in the last column shows the relative deficiency of our bound.\label{table:3}}\\[1cm]
\nopagebreak
\toprule
\nu & k & \utilde{j}'_{\nu,k} & \wtilde{j}'_{\nu,k} & \wtilde{\ell}'_{\nu,k} & \frac{\wtilde{j}'_{\nu,k}}{\wtilde{m}'_{\nu,k}}-1 \\\nopagebreak\midrule
\endfirsthead
\toprule
\nu & k & \utilde{j}'_{\nu,k} & \wtilde{j}'_{\nu,k} & \wtilde{\ell}'_{\nu,k} & \frac{\wtilde{j}'_{\nu,k}}{\wtilde{m}'_{\nu,k}}-1\\\nopagebreak\midrule
\endhead
 \multirow{9}{*}{$0$} & 2 & 3.82905543678956 & {\color{red}3.92699081698724} & 5.21994815723788 & \\\nopagebreak
 & 5 & 13.3236232806872 & {\color{red}13.3517687777566} & 17.8482723692095 & \\\nopagebreak
 & 10 & 29.0468218553286 & {\color{red}29.0597320457056} & 38.8629189129038 & \\\nopagebreak
 & 50 & 154.721014471994 & {\color{red}154.723438189297} & 206.941758158503 & \\\nopagebreak
 & 100 & 311.801868181958 & {\color{red}311.803070868787} & 417.036190247476 & \\\nopagebreak
 & 1000 & 3139.23633964380 & {\color{red}3139.23645909960} & 4198.72872262331 & \\\nopagebreak
 & 10000 & 31413.5703294702 & {\color{red}31413.5703414077} & 42015.6503274934 & \\\nopagebreak
 & 100000 & 314156.909163295 & {\color{red}314156.909164489} & 420184.866007090 & \\\nopagebreak
 & 500000 & 1570793.97060017 & {\color{red}1570793.97060041} & 2100936.93567796 & \\
\midrule
\multirow{10}{*}{$\frac{1}{2}$} & 1 &   & {\color{red}1.48584693535690} & 1.52657094660980 & \\\nopagebreak
 & 2 & 4.60247860677616 & {\color{red}4.68568659322651} & 5.84128905032192 & \\\nopagebreak
 & 5 & 14.1016664264204 & {\color{red}14.1283185385025} & 18.4716385900913 & \\\nopagebreak
 & 10 & 29.8283630264584 & {\color{red}29.8409412351793} & 39.4867451802621 & \\\nopagebreak
 & 50 & 155.505621002290 & {\color{red}155.508032535003} & 207.565903087165 & \\\nopagebreak
 & 100 & 312.586869472825 & {\color{red}312.588069144796} & 417.660372359456 & \\\nopagebreak
 & 1000 & 3140.02169802843 & {\color{red}3140.02181745436} & 4199.35293773371 & \\\nopagebreak
 & 10000 & 31414.3557236548 & {\color{red}31414.3557355921} & 42016.2745458797 & \\\nopagebreak
 & 100000 & 314157.694561061 & {\color{red}314157.694562255} & 420185.490225804 & \\\nopagebreak
 & 500000 & 1570794.75599825 & {\color{red}1570794.75599849} & 2100937.55989670 & \\
\midrule
\multirow{10}{*}{$1$} & 1 &   & 2.11506214188671 & {\color{red}2.11499393987502} & 3\times 10^{-5}\\\nopagebreak
 & 2 & 5.33004070160281 & {\color{red}5.40501381523099} & 6.45712461607948 & \\\nopagebreak
 & 5 & 14.8635367796392 & {\color{red}14.8889705581512} & 19.0933165595504 & \\\nopagebreak
 & 10 & 30.6019171965764 & {\color{red}30.6141946267452} & 40.1097896852094 & \\\nopagebreak
 & 50 & 156.288635764017 & {\color{red}156.291035345361} & 208.189900476501 & \\\nopagebreak
 & 100 & 313.371074961794 & {\color{red}313.372271647828} & 418.284481220551 & \\\nopagebreak
 & 1000 & 3140.80697683556 & {\color{red}3140.80709623164} & 4199.97714556513 & \\\nopagebreak
 & 10000 & 31415.1411098817 & {\color{red}31415.1411218186} & 42016.8987635385 & \\\nopagebreak
 & 100000 & 314158.479958031 & {\color{red}314158.479959224} & 420186.114444445 & \\\nopagebreak
 & 500000 & 1570795.54139618 & {\color{red}1570795.54139641} & 2100938.18411543 & \\
\midrule
\multirow{10}{*}{$5$} & 1 & 6.28826192531053 & 6.65074326003913 & {\color{red}6.54714840840174} & 2\times 10^{-2}\\\nopagebreak
 & 2 & 10.5188417338511 & {\color{red}10.5753107094334} & 11.2332989986105 & \\\nopagebreak
 & 5 & 20.5754812539651 & {\color{red}20.5957954835809} & 24.0099702812020 & \\\nopagebreak
 & 10 & 36.5607736396374 & {\color{red}36.5713809406576} & 45.0667270043476 & \\\nopagebreak
 & 50 & 162.498179300523 & {\color{red}162.500490908179} & 213.176588124942 & \\\nopagebreak
 & 100 & 319.616768720740 & {\color{red}319.617942511495} & 423.274719776078 & \\\nopagebreak
 & 1000 & 3147.08634964664 & {\color{red}3147.08646880501} & 4204.97054621798 & \\\nopagebreak
 & 10000 & 31421.4239132894 & {\color{red}31421.4239252240} & 42021.8924786219 & \\\nopagebreak
 & 100000 & 314164.763105141 & {\color{red}314164.763106335} & 420191.108190953 & \\\nopagebreak
 & 500000 & 1570801.82457384 & {\color{red}1570801.82457408} & 2100943.17786473 & \\
\midrule
\multirow{10}{*}{$10$} & 1 & 11.6398057744632 & 12.0162638908975 & {\color{red}11.8660736403637} & 1\times 10^{-2}\\\nopagebreak
 & 2 & 16.4468248709150 & {\color{red}16.4994512103844} & 16.9720812337165 & \\\nopagebreak
 & 5 & 27.1819924558840 & {\color{red}27.1998893619366} & 30.0345337375884 & \\\nopagebreak
 & 10 & 43.6067615870698 & {\color{red}43.6162303330504} & 51.1985421323618 & \\\nopagebreak
 & 50 & 170.135218972640 & {\color{red}170.137436706410} & 219.396825169749 & \\\nopagebreak
 & 100 & 327.357136966911 & {\color{red}327.358284407731} & 429.505960987985 & \\\nopagebreak
 & 1000 & 3154.92845527118 & {\color{red}3154.92857413496} & 4211.21164225011 & \\\nopagebreak
 & 10000 & 31429.2767018707 & {\color{red}31429.2767138023} & 42028.1345570098 & \\\nopagebreak
 & 100000 & 314172.616967415 & {\color{red}314172.616968609} & 420197.350367542 & \\\nopagebreak
 & 500000 & 1570809.67853160 & {\color{red}1570809.67853184} & 2100949.42005005 & \\
\midrule
\multirow{10}{*}{$50$} & 1 & 52.8187848373174 & 53.3254083348661 & {\color{red}53.0440529190390} & 5\times 10^{-3}\\\nopagebreak
 & 2 & 60.0249962636869 & {\color{red}60.0813194544867} & 60.2418339467916 & \\\nopagebreak
 & 5 & 74.3163203933786 & {\color{red}74.3314682956021} & 75.7386724516605 & \\\nopagebreak
 & 10 & 93.9425695282024 & {\color{red}93.9498510947580} & 98.4545059578280 & \\\nopagebreak
 & 50 & 227.750672604515 & {\color{red}227.752470650764} & 268.662409540753 & \\\nopagebreak
 & 100 & 387.108315157524 & {\color{red}387.109313313986} & 479.099715011984 & \\\nopagebreak
 & 1000 & 3217.38763707313 & {\color{red}3217.38775367737} & 4261.11426459507 & \\\nopagebreak
 & 10000 & 31492.0704532998 & {\color{red}31492.0704652076} & 42078.0685658962 & \\\nopagebreak
 & 100000 & 314235.445001727 & {\color{red}314235.445002920} & 420247.287518400 & \\\nopagebreak
 & 500000 & 1570872.50962077 & {\color{red}1570872.50962101} & 2100999.35748022 & \\
\midrule
\multirow{10}{*}{$100$} & 1 & 103.552648845966 & 104.160029957632 & {\color{red}103.803049301258} & 3\times 10^{-3}\\\nopagebreak
 & 2 & 112.385101849666 & {\color{red}112.448430986771} & 112.528002221252 & \\\nopagebreak
 & 5 & 129.358651207261 & {\color{red}129.374174936818} & 130.325058224489 & \\\nopagebreak
 & 10 & 151.813259002596 & {\color{red}151.820223347072} & 155.042104955214 & \\\nopagebreak
 & 50 & 294.664327855513 & {\color{red}294.665911746948} & 329.158665247029 & \\\nopagebreak
 & 100 & 457.918838435071 & {\color{red}457.919731860974} & 540.487440358064 & \\\nopagebreak
 & 1000 & 3294.79831773381 & {\color{red}3294.79843173622} & 4323.42741287888 & \\\nopagebreak
 & 10000 & 31570.4915863354 & {\color{red}31570.4915982138} & 42140.4795335826 & \\\nopagebreak
 & 100000 & 314313.972888315 & {\color{red}314313.972889508} & 420309.708302353 & \\\nopagebreak
 & 500000 & 1570951.04705006 & {\color{red}1570951.04705030} & 2101061.77913701 & \\
\midrule
\multirow{10}{*}{$1000$} & 1 & 1007.65497214255 & 1008.87692662131 & {\color{red}1008.10741918263} & 8\times 10^{-4}\\\nopagebreak
 & 2 & 1025.97327758496 & 1026.08853535440 & {\color{red}1026.00990668415} & 8\times 10^{-5}\\\nopagebreak
 & 5 & 1059.53767346269 & {\color{red}1059.56146481700} & 1059.76655932628 & \\\nopagebreak
 & 10 & 1101.19242597862 & {\color{red}1101.20147226161} & 1102.01144720199 & \\\nopagebreak
 & 50 & 1326.56778093642 & {\color{red}1326.56910593182} & 1338.72297124053 & \\\nopagebreak
 & 100 & 1546.19139906085 & {\color{red}1546.19204915940} & 1582.10105006399 & \\\nopagebreak
 & 1000 & 4600.92489073707 & {\color{red}4600.92497958265} & 5433.32953880577 & \\\nopagebreak
 & 10000 & 32969.1998239474 & {\color{red}32969.1998353403} & 43262.6382055024 & \\\nopagebreak
 & 100000 & 315726.121837736 & {\color{red}315726.121838924} & 421433.158077796 & \\\nopagebreak
 & 500000 & 1572364.44893452 & {\color{red}1572364.44893476} & 2102185.34408437 & \\
\midrule
\multirow{10}{*}{$10000$} & 1 & 10016.4918376026 & 10019.0850161927 & {\color{red}10017.4305732935} & 2\times 10^{-4}\\\nopagebreak
 & 2 & 10055.6275626539 & 10055.8663257080 & {\color{red}10055.6423608815} & 2\times 10^{-5}\\\nopagebreak
 & 5 & 10126.5382295259 & {\color{red}10126.5852535008} & 10126.5888462274 & \\\nopagebreak
 & 10 & 10213.1216137128 & {\color{red}10213.1385542345} & 10213.3033978649 & \\\nopagebreak
 & 50 & 10658.2856359516 & {\color{red}10658.2876108152} & 10661.1849966748 & \\\nopagebreak
 & 100 & 11062.1238234066 & {\color{red}11062.1246535611} & 11071.3066333350 & \\\nopagebreak
 & 1000 & 15489.7059732163 & {\color{red}15489.7060377999} & 15852.5040437567 & \\\nopagebreak
 & 10000 & 46030.9746376957 & {\color{red}46030.9746465753} & 54361.7461824340 & \\\nopagebreak
 & 100000 & 329713.213863324 & {\color{red}329713.213864463} & 432654.745761900 & \\\nopagebreak
 & 500000 & 1586470.41726056 & {\color{red}1586470.41726079} & 2113418.40311191 & \\
\midrule
\multirow{10}{*}{$100000$} & 1 & 100035.530371101 & 100041.098972602 & {\color{red}100037.536956311} & 4\times 10^{-5}\\\nopagebreak
 & 2 & 100119.692630670 & 100120.202625833 & {\color{red}100119.708964190} & 5\times 10^{-6}\\\nopagebreak
 & 5 & 100271.813981416 & 100271.913355254 & {\color{red}100271.825350454} & 9\times 10^{-7}\\\nopagebreak
 & 10 & 100456.883755739 & {\color{red}100456.919097553} & 100456.923246477 & \\\nopagebreak
 & 50 & 101397.105043612 & {\color{red}101397.108914265} & 101397.740786006 & \\\nopagebreak
 & 100 & 102234.528924186 & {\color{red}102234.530469187} & 102236.567493622 & \\\nopagebreak
 & 1000 & 110670.793673338 & {\color{red}110670.793755672} & 110763.646675313 & \\\nopagebreak
 & 10000 & 154924.826290601 & {\color{red}154924.826297055} & 158556.516370227 & \\\nopagebreak
 & 100000 & 460331.471036162 & {\color{red}460331.471037050} & 543645.912187864 & \\\nopagebreak
 & 500000 & 1724974.19659187 & {\color{red}1724974.19659209} & 2225491.96489915 & \\
\midrule
\multirow{10}{*}{$500000$} & 1 & 500060.756008294 & 500070.272554347 & {\color{red}500064.182433034} & 1\times 10^{-5}\\\nopagebreak
 & 2 & 500204.624142480 & 500205.494863009 & {\color{red}500204.649221081} & 2\times 10^{-6}\\\nopagebreak
 & 5 & 500464.546806439 & 500464.716136887 & {\color{red}500464.551147341} & 3\times 10^{-7}\\\nopagebreak
 & 10 & 500780.557830912 & 500780.617909442 & {\color{red}500780.571426180} & 9\times 10^{-8}\\\nopagebreak
 & 50 & 502382.483279350 & {\color{red}502382.489781183} & 502382.701278916 & \\\nopagebreak
 & 100 & 503804.363624797 & {\color{red}503804.366193259} & 503805.063905356 & \\\nopagebreak
 & 1000 & 517885.834593690 & {\color{red}517885.834718306} & 517918.417812500 & \\\nopagebreak
 & 10000 & 586227.229833856 & {\color{red}586227.229841076} & 587665.712270748 & \\\nopagebreak
 & 100000 & 967170.501022458 & {\color{red}967170.501023138} & 1017601.68651123 & \\\nopagebreak
 & 500000 & 2301667.01053780 & {\color{red}2301667.01053798} & 2718242.20552639 & \\
\bottomrule
\end{longtable}}

\begin{figure}[ht]
\centering
\includegraphics{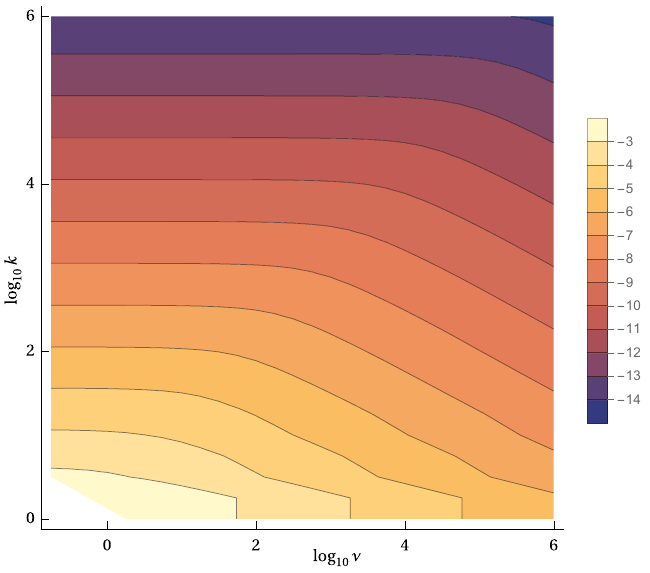}
\caption{The contour plot of $\log_{10}\left(\frac{\wtilde{j}'_{\nu,k} }{\utilde{j}'_{\nu,k}}-1\right)$, which estimates from above the maximum possible order of the relative error of our approximations of $j'_{\nu,k}$.\label{fig:errjprime}}
\end{figure}

\end{document}